\documentclass[runningheads]{llncs}
\usepackage{latexsym}
\def\ra{\rightarrow}
\def\lra{\leftrightarrow}
\def\llra{\longleftrightarrow}
\def\be{\begin{equation}}
\def\ee{\end{equation}}
\def\bea{\begin{eqnarray}}
\def\eea{\end{eqnarray}}
\def\A{{\sf A}} %unfortunately plain A is an abelian surface for now
\def\tr{{\rm tr}}
\def\Nm{{\rm N}}
\def\H{{\cal H}}
\def\O{{\cal O}}
\def\F{{\bf F}}
\def\Z{{\bf Z}}
\def\C{{\bf C}}
\def\Q{{\bf Q}}
\def\R{{\bf R}}
\def\Pr{{\bf P}}
\def\Aut{{\rm Aut}}
\def\PSL{{\rm PSL}}
\def\PGL{{\rm PGL}}
\def\calX{{\cal X}}
\def\calXo{\calX_0}
\def\Xone{{\cal X}(1)}
\def\XX{{\cal X}^*}
\def\XXo{{\XX_0}}
\def\ww{{\sf w}}
\def\geqs{\geq}
\def\leqs{\leq}
\begin{document}

\title{Shimura Curve Computations}
\author{Noam D.~Elkies}
\institute{Harvard University}
\maketitle

\begin{abstract}
We give some methods for computing equations for certain Shimura
curves, natural maps between them, and special points on them.
We then illustrate these methods by working out several examples
in varying degrees of detail.  For instance, we compute coordinates
for all the rational CM points on the curves $\XX(1)$ associated
with the quaternion algebras over~$\Q$ ramified at $\{2,3\}$,
$\{2,5\}$, $\{2,7\}$, and $\{3,5\}$.  We conclude with a list of
open questions that may point the way to further computational
investigation of these curves.
\end{abstract}
\section{Introduction}
\subsection{Why and how to compute with Shimura curves}

The classical modular curves, associated to congruence subgroups
of $\PSL_2(\Q)$, have long held and repaid the interest of number
theorists working theoretically as well as computationally.
In the fundamental paper~\cite{S2} Shimura defined curves associated
with other quaternion algebras other over totally real number
fields in the same way that the classical curves are associated
with the algebra $M_2(\Q)$ of $2\times 2$ matrices over~$\Q$.  These
Shimura curves are now recognized as close analogues of the classical
modular curves: almost every result involving the classical curves
generalizes with some more work to Shimura curves, and indeed Shimura
curves figure alongside classical ones in a key step in the recent
proof of Fermat's ``last theorem''~\cite{Ribet}.

But computational work on Shimura curves lags far behind the extensive
effort devoted to the classical modular curves.  The 19th century
pioneers investigated some arithmetic quotients of the upper half plane
which we now recognize as Shimura curves (see for instance \cite{F1,F2})
with the same enthusiasm that they applied to the $\PSL_2(\Q)$ curves.
But further inroads proved much harder for Shimura curves than for
their classical counterparts.  The $\PSL_2(\Q)$ curves parametrize
elliptic curves with some extra structure; the general elliptic
curve has a simple explicit formula which lets one directly write
down the first few modular curves and maps between them.  (For
instance, this is how Tate obtained the equations for the first
few curves X$_1(N)$ parametrizing elliptic curves with an $N$\/-torsion
point; see for instance \cite[pp.145--148]{Knapp}.)  Shimura showed
that curves associated with other quaternion algebras also
parametrize geometric objects, but considerably more complicated ones
(abelian varieties with quaternionic endomorphisms); even in the first
few cases beyond $M_2(\Q)$, explicit formulas for these objects were
obtained only recently~\cite{HM}, and using such formulas to get at
the Shimura curves seems a most daunting task.  Moreover, most modern
computations with modular curves (e.g.~\cite{Cr,Modcomp}) sidestep the
elliptic interpretation and instead rely heavily on $q$-expansions,
i.e.\ on the curves' cusps.  But arithmetic subgroups of $\PSL_2(\R)$
other than those in $\PSL_2(\Q)$ contain no parabolic elements,
so their Shimura curves have no cusps, and thus any method that
requires $q$-expansions must fail.  

But while Shimura curves pose harder computational problems than
classical modular curves, efficient solutions to these problems
promise great benefits.  These curves tempt the computational
number theorist not just because, like challenging mountainpeaks,
``they're there''$\!$, but because of their remarkable properties,
direct applications, and potential for suggesting new ideas for
theoretical research.  Some Shimura curves and natural maps between
them provide some of the most interesting examples in the geometry
of curves of low genus; for instance each of the five curves of
genus $g\in[2,14]$ that attains the Hurwitz bound $84(g-1)$ on the
number of automorphisms of a curve in characteristic~zero is a
Shimura curve.  Shimura curves, like classical and Drinfeld modular
curves, reduce to curves over the finite field $\F_{q^2}$ of $q^2$
elements that attain the Drinfeld-\mbox{\kern-.16em}Vl\u{a}du\c{t}
upper bound $(q-1+o(1)) g$ on the number of points
of a curve of genus $g\ra\infty$ over that field~\cite{I3}.
Moreover, while all three flavors of modular curves include
towers that can be given by explicit formulas and thus used
to construct good error-correcting codes~\cite{Goppa1,Goppa2,TVZ},
only the Shimura curves, precisely because of their lack of cusps,
can give rise to totally unramified towers, which should simplify the
computation of the codes; we gave formulas for several such towers
in~\cite{Modtower}.  Finally, the theory of modular curves indicates
that CM (complex multiplication) points on Shimura curves, elliptic
curves covered by them, and modular forms on them have number-theoretic
significance.  The ability to efficiently compute such objects should
suggest new theoretical results and conjectures concerning the
arithmetic of Shimura curves.  For instance, the computations of
CM~points reported in this paper should suggest factorization
formulas for the difference between the coordinates of two
such points analogous to those of Gross and Zagier~\cite{GZ}
for $j$-invariants of elliptic curves, much as the computation of 
CM values of the Weber modular functions suggested the formulas
of~\cite{YZ}.  Also, as in \cite{GS}, rational CM points on rational
Shimura curves with only three elliptic points (i.e.\ coming from
arithmetic triangle groups $G_{p,q,r}$) yield identities $A+B=C$\/
in coprime integers $A,B,C$\/ with many repeated factors; we list the
factorizations here, though we found no example in which
$A,B,C$\/ are perfect $p,q,r$-th powers, nor any new near-record
ABC ratios.  Finally, CM computations on Shimura curves may also make
possible new Heegner-point constructions as in~\cite{Heegcomp}.

So how do we carry out these computations?  In a few cases (listed
in~\cite{JL}), the extensive arithmetic theory of Shimura curves has
been used to obtain explicit equations, deducing from the curves'
$p$-adic uniformizations Diophantine conditions on the coefficients
of their equations stringent enough to determine them uniquely.
But we are interested, not only in the equations, but in modular
covers and maps between Shimura curves associated to the same
quaternion algebra, and in CM~points on those curves.  The arithmetic
methods may be able to provide this information, but so far no such
computation seems to have been done.  Our approach relies mostly
on the uniformization of these curves {\em qua} Riemann surfaces
by the hyperbolic plane, and uses almost no arithmetic.  This
approach is not fully satisfactory either; for instance it probably
cannot be used in practice to exhibit all natural maps between Shimura
curves of low genus.  But it will provide equations for at least a
hundred or so curves and maps not previously accessible, which include
some of the most striking examples and should provide more than enough
data to suggest further theoretical and computational work.

When a Shimura curve~$C$\/ comes from an arithmetic subgroup
of~$\PSL_2(\R)$ contained in a triangle group $G_{p,q,r}$, the curve
$\H/G_{p,q,r}$ has genus~0, and $C$\/ is a cover of that curve branched
only above three points, so may be determined from the ramification
data.  (We noted in~\cite[p.48]{Modcomp} that this method was
available also for classical modular curves comings from subgroups
of $\PSL_2(\Z)\cong G_{2,3,\infty}$, though there better methods
are available thanks to the cusp.  Subgroups of~$\PSL_2(\R)$
commensurate with\footnote{
  Recall that two subgroups $H,K$\/ of a group~$G$\/ are said
  to be {\em commensurate} if $H\cap K$\/ is a subgroup of finite
  index in both~$H$\/ and~$K$.
  }
but not contained in $G_{p,q,r}$ may be handled
similarly via the common subgroup of finite index.)  The
identification of $\H/G_{p,q,r}$ with $\Pr^1$ is then
given by a quotient of hypergeometric functions on $\Pr^1$,
which for instance lets us compute the $\Pr^1$ coordinate of
any CM point on~$C$\/ as a complex number to high precision and
thus recognize it at least putatively as an algebraic number.

Now it is known~\cite{Takeuchi} that only nineteen commensurability
classes of arithmetic subgroups of $\PSL_2(\R)$ contain a triangle
group.  These include some of the most interesting examples --- for
instance, congruence subgroups of arithmetic triangle groups account
for several of the sporadic ``arithmetically exceptional functions''
(rational functions $f(X)\in\Q(X)$ which permute $\Pr^1(\F_{\!p})$
for infinitely many primes~$p$) of~\cite{Muller}; but
an approach that could only deal with those nineteen classes would
be limited indeed.  When there are more than three elliptic points,
a new difficulty arises: even if $C=\cal H/G$\/ still has genus~0,
we must first determine the relative locations of the elliptic points,
and to locate other CM points we must replace the hypergeometric
functions to solutions of more general ``Schwarzian differential
equations'' in the sense of~\cite{I1}.  We do both by in effect
using nontrivial elements of the ``commensurator'' of the group
$G\in\PSL_2(\R)$, i.e.\ transformations in $\PSL_2(\R)$ which
do not normalize~$G$\/ but conjugate $G$\/ to a group commensurable
with~$G$.  Ihara had already used these commensurators in~\cite{I1}
theoretically to prove that both $C$\/ and its Schwarzian equation
are defined over a number field, but this method has apparently
not been actually used to compute such equations until now.

\subsection{Overview of the paper}

We begin with a review of the necessary definitions and facts
on quaternion algebras and Shimura curves, drawn mostly from
\cite{S2} and~\cite{V}.  We then give extended computational accounts
of Shimura curves and their supersingular and rational CM points
for the two simplest indefinite quaternion algebras over~$\Q$
beyond the classical case of the matrix algebra $M_2(\Q)$, namely
the quaternion algebras ramified at $\{2,3\}$ and $\{2,5\}$.
In the final section we more briefly treat some other examples
which illustrate features of our methods that do not arise in the
$\{2,3\}$ and $\{2,5\}$ cases, and conclude with some open questions
suggested by our computations that may point the way to further
computational investigation on these curves.

\subsection{Acknowledgements}  Many thanks to B.H.~Gross for
introducing me to Shimura curves and for many enlightening conversations
and clarifications on this fascinating topic.  Thanks also to Serre
for a beautiful course that introduced me to three-point covers of
$\Pr^1$ among other things (\cite{Serre}, see also \cite{Mat}); to
Ihara for alerting me to his work~\cite{I1,I2} on supersingular points
on Shimura curves and their relation with the curves' uniformization
by the upper half-plane; and to C.~McMullen for discussions of the
uniformization of quotients of $\H$ by general co-compact
discrete subgroups of $\PSL_2(\R)$.  A.~Adler provided several
references to the 19th-century literature, and C.~Doran informed
me of~\cite{HM}.  Finally, I thank B.~Poonen for reading and
commenting on a draft of this paper, leading to considerable
improvements of exposition in several places.

The numerical and symbolic computations reported here were 
carried out using the {\sc gp/pari} and {\sc macsyma} packages, except
for (\ref{14:5}), for which I thank Peter M\"uller as noted there.

This work was made possible in part by funding from the
David and Lucile Packard Foundation.

\section{Review of quaternion algebras over $\Q$
and their Shimura curves}

\subsection{Quaternion algebras over~$\Q$; the arithmetic groups
$\Gamma(1)$ and $\Gamma^*(1)$}

Let $K$\/ be a field of characteristic zero; for our purposes
$K$\/ will always be a number field or, rarely, its localization,
and usually the number field will be~$\Q$.  A {\em quaternion algebra
over~$K$}\/ is a simple associative algebra~$\A$ with unit, containing~$K$,
such that $K$\/ is the center of~$\A$ and $\dim_K\A=4$.
Such an algebra has a {\em conjugation} $a\lra\bar a$,
which is a $K$-linear anti-involution (i.e.\ $\bar{\bar a}=a$
and $\overline{a_1 a_2}=\bar a_2\bar a_1$ hold identically in~$\A$)
such that $a=\bar a \Leftrightarrow a\in K$.
The {\em trace} and {\em norm} are the additive and multiplicative
maps from~$\A$ to~$K$\/ defined by
\be
\tr(a) = a + \bar a,\qquad
\Nm(a) = a \bar a = \bar a a;
\label{t,N}
\ee
every $a\in\A$ satisfies its characteristic equation
\be
a^2 - (\tr(a)) a + \Nm(a) = 0.
\label{chareq}
\ee
The most familiar example of a quaternion algebra is $M_2(K)$, the algebra
of $2\times 2$ matrices over~$K$, and if $K$\/ is algebraically closed
then $M_2(K)$ is the only quaternion algebra over~$K$\/ up to isomorphism.
The other well-known example is the algebra of Hamilton quaternions
over~$\R$.  In $M_2(K)$ the trace is the usual trace of a square matrix,
so the conjugate of $a\in M_2(K)$ is $\tr(a)I_{2\times 2} - a$,
and the norm is just the determinant.  Any quaternion algebra with
zero divisors is isomorphic with~$M_2(K)$.  An equivalent criterion
is that the algebra contain a nonzero element whose norm and trace
both vanish.  Now the trace-zero elements constitute a $K$-subspace
of~$\A$ of dimension~3, on which the norm is a homogeneous quadric;
so the criterion states that $\A\cong M_2(K)$ if and only if
that quadric has nonzero $K$-rational points.  The Hamilton quaternions
have basis $1,i,j,k$\/ satisfying the familiar relations
\be
i^2 = j^2 = k^2 = 1, \quad
ij = -ji = k, \quad
jk = -kj = i, \quad
ki = -ik = j;
\label{Hamilton}
\ee
the conjugates of $1,i,j,k$\/ are $1,-i,-j,-k$, so a Hamilton
quaternion $\alpha_1 + \alpha_2 i + \alpha_3 j + \alpha_4 k$\/
has trace $2\alpha_1$ and norm 
$\alpha_1^2 + \alpha_2^2 + \alpha_3^2 + \alpha_4^2$.  Thus the
Hamilton quaternions over~$K$\/ are isomorphic with $M_2(K)$
if and only if $-1$ is a sum of two squares in~$K$.

In fact it is known that if $K=\R$ then every quaternion algebra
over~$K$\/ is isomorphic with either $M_2(\R)$ or the Hamilton
quaternions.  In general if $K$\/ is any local field of
characteristic~zero then there is up to isomorphism exactly one
quaternion algebra over~$K$\/ other than $M_2(K)$ --- with the
exception of the field of complex numbers, which being algebraically
closed admits no quaternion algebras other than~$M_2(\C)$.  If $\A$ is
a quaternion algebra over a number field~$K$\/ then a finite or
infinite place~$v$ of~$K$\/ is said to be {\em ramified}\/ in~$\A$
if $\A\otimes K_v$ is not isomorphic with $M_2(K_v)$.  There can
only be a finite number of ramified places, because a nondegenerate
quadric over~$K$\/ has nontrivial local zeros at all but finitely many
places of~$K$.  A less trivial result (the case $K=\Q$ is equivalent
to Quadratic Reciprocity) is that the number of ramified places is
always even, and to each finite set of places~$\Sigma$ of even
cardinality containing no complex places there corresponds a unique
(again up to isomorphism) quaternion algebra over~$K$\/ ramified
at those places and no others.  In particular an everywhere unramified
quaternion algebra over~$K$\/ must be isomorphic with $M_2(K)$.

An {\em order} in a quaternion algebra over a number field (or a
non-Archimedean local field)~$K$\/ is a subring containing the ring
$O_K$ of $K$-integers and having rank~4 over~$O_K$.  For instance
$M_2(O_K)$ and $O_K[i,j]$ are orders in the matrix and quaternion
algebras over~$K$.  Any order is contained in at least one
{\em maximal}\/ order, that is, in an order not properly contained
in any other.  Examples of maximal orders are $M_2(O_K)\in M_2(K)$ and
the Hurwitz order $\Z[1,i,j,(1+i+j+k)/2]$ in the Hamilton quaternions
over~$\Q$.  It is known that if $K$\/ has at least one Archimedean
place at which $\A$ is {\em not}\/ isomorphic with the Hamilton
quaternions then all maximal orders are conjugate in~$\A$.

Now let\footnote{
  Most of our examples, including the two that will occupy us in
  the next two sections, involve quaternion algebras over~$\Q$.
  In~\cite{S2} Shimura associated modular curves to a quaternion
  algebra over any totally real number field~$K$\/ for which the
  algebra is ramified at all but one of the infinite places of~$K$.
  Since the special case $K=\Q$ accounts for most of our computations,
  and is somewhat easier to describe, we limit our discussion to
  quaternion algebras over~$\Q$ from here until section 5.3.  At
  that point we briefly describe the situation for arbitrary~$K$\/
  before working out a couple of examples with $[K:\Q]>$1.
  }
$K=\Q$.  A quaternion algebra $\A/\Q$ is called {\em definite}
or {\em indefinite} according as $\A\otimes\R$ is isomorphic with the
Hamilton quaternions or $M_2(\R)$, i.e.\ according as the infinite
place is ramified or unramified in~$\A$. [These names allude to the
norm form on the trace-zero subspace of~$\A$, which is definite
in the former case, indefinite in the latter.]  We shall be concerned
only with the indefinite case.  Then $\Sigma$ consists of an even number
of finite primes.  Fix such a $\Sigma$ and the corresponding quaternion
algebra~$\A$.  Let $\O$ be a maximal order in~$\A$; since $\A$ is
indefinite, all its maximal orders are conjugate, so choosing a
different maximal order would not materially affect the constructions
in the sequel.

Let $\O^*_1$ be the group of units of norm~1 in~$\O$.
We then define the following arithmetic subgroups
of $\A^*/\Q^*$:
\bea
\Gamma(1) & := & \O^*_1/\{\pm1\},
\label{G}
\\
\Gamma^*(1) & := &
\{ [a] \in \A^*/\Q^*: a \O = \O a, \ \Nm(a) > 0 \}.
\label{G*}
\eea
[In other words $\Gamma^*(1)$ is the normalizer of $\Gamma(1)$ in the
positive-norm subgroup of $\A^*/\Q^*$.  Takeuchi~\cite{Takeuchi} calls
these groups $\Gamma^{(1)} (\A,\O_1)$ and $\Gamma^{(*)} (\A,\O_1)$;
we use $\Gamma(1)$ to emphasize the analogy with the classical
case of $\PSL_2(\Z)$, which makes $\Gamma^*(1)$ a natural
adaptation of Takeuchi's notation.  Vign\'eras \cite[p.~121ff.]{V}
calls the same groups $\Gamma$ and~$G$, citing~\cite{Mi} for the
structure of their quotient.]
As noted, $\Gamma(1)$ is a normal subgroup of $\Gamma^*(1)$.
In fact $\Gamma^*(1)$ consists of the classes mod~$\Q^*$
of elements of~$\O$ whose norm is $\prod_{p\in \Sigma'} p$
for some (possibly empty) subset $\Sigma'\subseteq \Sigma$,
and $\Gamma^*(1)/\Gamma(1)$ is an elementary abelian
2-group with $\#\Sigma$ generators.

\subsection{The Shimura modular curves $\Xone$ and $\XX(1)$}

The group $\Gamma(1)$, and thus any other group commensurable
with it such as $\Gamma^*(1)$, is a discrete subgroup of
$(A\otimes\R)^*_+/\R^*$ (the subscript ``$+$'' indicating
positive norm), with compact quotient unless $\Sigma=\emptyset$, and
of finite covolume even in that case.  Since $A\otimes\R\cong M_2(\R)$,
the group $(A\otimes\R)^*_+/\R^*$ is isomorphic with $\PSL_2(\R)$ and
thus with $\Aut(\H)$, the group of automorphisms of the hyperbolic
upper half plane
\be
\H := \{ z\in\C: {\rm Im}(z) > 0 \}.
\label{Hdef}
\ee
Explicitly, a unimodular matrix $\pm({a\;b\atop c\;d})$ acts on~$\H$
via the fractional linear transformation $z \mapsto (az+b) / (cz+d)$.
We may define the Shimura curves $\Xone$ and $\XX(1)$ {\em qua} compact
Riemann surfaces by
\be
\Xone := \H/\Gamma(1),\qquad
\XX(1) := \H/\Gamma^*(1).
\label{Xdef}
\ee
[More precisely, the Riemann surfaces are given by (\ref{Xdef}) unless
$\Sigma=\emptyset$, in which case the quotient only becomes compact
upon adjoining a cusp.]  The hyperbolic area of these quotients
of~$\H$ is given by the special case $k=\Q$ of a formula of
Shimizu~\cite[Appendix]{S1}, quoted in \cite[p.207]{Takeuchi}.
Using the normalization
$\pi^{-1}\int\!\!\!\int dx\,dy/y^2$ for the hyperbolic area
(with $z=x+iy$; this normalization gives an ideal triangle unit area),
that formula is
\be
{\rm Area}(\Xone) = \frac16 \prod_{p\in\Sigma} (p-1),
\label{Area1}
\ee
from which
\be
{\rm Area}(\XX(1)) = \frac1{[\Gamma^*(1):\Gamma(1)]} {\rm Area}(\Xone)
= \frac16 \prod_{p\in\Sigma} \frac{p-1}{2}.
\label{Area*}
\ee

It is known (see for instance Ch.IV:2,3 of~\cite{V} for the following
facts) that, for any discrete subgroup $\Gamma\subset\PSL_2(\R)$
of finite covolume, the genus of $\H/\Gamma$ is determined by its
area together with with information on elements of finite order
in~$\Gamma$.  All finite subgroups of~$\Gamma$ are cyclic, and there
are finitely many such subgroups up to conjugation in~$\Gamma$.
There are finitely many points $P_j$ of $\H/\Gamma$ with nontrivial
stabilizer, and the stabilizers are the maximal nontrivial finite
subgroups of~$\Gamma$ modulo conjugation in~$\Gamma$.  If the
order of the stabilizer of~$P_j$ is $e_j$ then $P_j$ is said to
be an ``elliptic point of order~$e_j$''$\!$.  Then if $\H/\Gamma$
is compact then its genus $g=g(\H/\Gamma)$ is given by
\be
2 g - 2 = {\rm Area}(\H/\Gamma) - \sum_j (1 - \frac1{e_j}).
\label{genus}
\ee
Moreover $\Gamma$ has a presentation
\be
\Gamma = \langle \alpha_1,\ldots,\alpha_g,\beta_1,\ldots,\beta_g,s_j |
s_j^{e_j} = 1,\
\prod_j s_j \prod_{i=1}^g [\alpha_i,\beta_i] = 1\rangle,
\label{G_present}
\ee
in which $s_j$ generates the stabilizer of a preimage of~$P_j$
in~$\H$ and rotates a neighborhood of that preimage by an angle
$2\pi/e_j$ (i.e.\ has derivative $e^{2\pi i/e_j}$ at its fixed point),
and $[\alpha,\beta]$ is the commutator
$\alpha\beta\alpha^{-1}\beta^{-1}$.
[This group is sometimes called $(g;e_1,\ldots,e_g)$.]
If $\H/\Gamma$ is not compact then we must subtract the number of
cusps from the right-hand side of~(\ref{genus}) and include a
generator~$s_j$ of~$\Gamma$ of infinite order for each cusp, namely
a generator of the infinite cyclic stabilizer of the cusp.  This
generator is a ``parabolic element'' of~$\PSL_2(\R)$, i.e.\
a fractional linear transformation with a single fixed point;
there are two conjugacy classes of such elements in $\PSL_2(\R)$,
and $s_j$ will be in the class of $z\mapsto z+1$.  We assign
$e_j=\infty$ to a cusp.  For both finite and infinite~$e_j$,
the trace and determinant of~$s_j$ are related by
\be
{\rm Tr}^2(s_j) = 4 \cos^2 \frac{\pi}{e_j} \det(s_j).
\label{tr-det}
\ee
Since we are working in quaternion algebras over~$\Q$, this means
that $e_j\in\{2,3,4,6,\infty\}$, and only $2,3,\infty$ are possible
if $\Gamma\subseteq\Gamma(1)$.  Moreover $e_j=\infty$ occurs only
in the classical case $\Sigma=\emptyset$.

We shall need to numerically compute for several such~$\Gamma$
the identification of $\H/\Gamma$ with an algebraic curve $X/\C$,
i.e.\ to compute the coordinates on~$X$\/ of a point corresponding
to (the $\Gamma$-orbit of) a given $z\in\H$, or inversely to obtain
$z$ corresponding to a point with given coordinates.  In fact the
two directions are essentially equivalent, because if we can
efficiently compute an isomorphism between two Riemann surfaces
then we can compute its inverse almost as easily.  For classical
modular curves one usually uses $q$-expansions to go from~$z$ to
rational coordinates; but this method is not available for our
groups~$\Gamma$, which have no parabolic ($e_j=\infty$) generator.
We can, however, still go in the opposite direction, computing
the map from~$X$\/ to $\H/\Gamma$ by solving differential equations
on~$X$.  The key is that while the function~$z$ on~$X$\/ is not well
defined due to the $\Gamma$ ambiguity, its Schwarzian derivative is.
In local coordinates the {\em Schwarzian derivative} of a nonconstant
function $z=z(\zeta)$ is the meromorphic function defined by
\be
S_\zeta(z) :=
-4 z^{-1} {z'}^{1/2} \frac{d^2}{d\zeta^2} \frac{z}{{z'}^{1/2}}
= \frac{2 z' z''' - 3 {z''}^2} {{z'}^2}.
\label{Schwarz}
\ee
This vanishes if and only if $z$ is a fractional linear transformation
of~$\zeta$.  Moreover it satisfies a nice ``chain rule'': if $\zeta$
is in turn a function of~$\eta$ then
\be
S_\eta(z) =
\left( \frac{d\zeta}{d\eta} \right)^2 S_\zeta(z) + S_\eta(\zeta).
\label{Schwarz:chain}
\ee
Thus if we choose a coordinate~$\zeta$ on~$X$\/ then $S_\zeta(z)$
is the same for each lift of~$z$ from~$\H/\Gamma$ to~$\H$, and thus
gives a well-defined function on the complement in~$X$\/ of the
elliptic points; changing the coordinate from~$\zeta$ to~$\eta$
multiples this function by $(d\zeta/d\eta)^2$ and adds a term
$S_\eta(\zeta)$ that vanishes if $\zeta$ is a fractional linear
transformation of~$\eta$.  In particular if $X$\/ has genus~0 and
we choose only rational coordinates (i.e.\ $\eta,\zeta$ are rational
functions of degree~1) then these terms $S_\eta(\zeta)$ always
vanish and $S_\zeta(z)\,d\zeta^2$ is a well-defined quadratic
differential~$\sigma$ on~$X$.  Near an elliptic point~$\zeta_0$ of
index~$e_j$, the function~$z$ has a branch point such that
$(z-z_0)/(z-\bar z_0)$ is $(\zeta-\zeta_0)^{1/e_j}$ times an analytic
function; for such~$z$ the Schwarzian derivative is still well-defined
in a neighborhood of~$\zeta_0$ but has a double pole there with leading
term $(1-e_j^{-2})/(\zeta-\zeta_0)^2$ [or $(1-e_j^{-2})/\zeta^2$
if $\zeta_0=\infty$ --- note that this too has a double pole
when multiplied by $d\zeta^2$].  So $\sigma = S_\zeta(z)\,d\zeta^2$
is a rational quadratic differential on~$X$, regular except for
double poles of known residue at the elliptic points, and independent
of the choice of rational coordinate when $X$\/ has genus~0.
Knowing $\sigma$ we may recover $z$ from the differential equation
\be
S_\zeta(z) = \sigma/d\zeta^2,
\label{dfq3}
\ee
which determines $z$ up to a fractional linear transformation
over~$\C$, and can then remove the ambiguity if we know at least
three values of~$z$ (e.g.\ at elliptic points, which are fixed
points of known elements of~$\Gamma$).

Because $S_\zeta(z)$ is invariant under fractional linear
transformations of~$z$, the third-order nonlinear differential
equation~(\ref{dfq3}) can be linearized as follows (see e.g.\
\cite[\S1--5]{I1}).  Let $(f_1,f_2)$ be a basis for the solutions of
the linear second-order equation
\be
f'' = a f' + b f
\label{dfq2}
\ee
for some functions $a(\zeta),b(\zeta)$.  Then
$z:=f_1/f_2$ is determined up to fractional linear transformation,
whence $S_\zeta(z)$ depends only on $a,b$ and not the
choice of basis.  In fact we find, using either of the equivalent
definitions in~(\ref{Schwarz}), that
\be
S_\zeta(f_1/f_2) = 2\frac{da}{d\zeta} -  a^2 - 4b.
\label{schw.cond}
\ee
Thus if $a$ is any rational function and
$b = -\sigma/4 d\zeta^2 + a'/2 - a^2/4$
then the solutions of~(\ref{dfq3}), and thus a map
from~$X$\/ to $\H/\Gamma$, are ratios of linearly independent
pairs of solutions of~(\ref{dfq2}).  In the terminology of~\cite{I1},
(\ref{dfq2}) is then a {\em Schwarzian equation} for $\H/\Gamma$.
We shall always choose $a$ so that $a\,d\zeta$ has
at most simple poles at the elliptic points and no other poles;
the Schwarzian equation then has regular singularities at the
elliptic points and no other singularities.
The most familiar example is the case that $\Gamma$ is a triangle
group, i.e.\ $X$\/ has genus~0 and three elliptic points (if $g=0$
there must be at least three elliptic points by~(\ref{genus})).  In
that case $\sigma$ is completely determined by its poles and residues:
if two different $\sigma$'s were possible, their difference would
be a nonzero quadratic differential on~$\Pr^1$ with at most three
simple poles, which is impossible.  If we choose the coordinate
on~$X$\/ that puts the elliptic points at $0,1,\infty$, and
require that $a$ be chosen of the form $a = C_0/\zeta + C_1/(\zeta-1)$
so that $b$ has only simple poles at~$0,1$, then there are four
choices for $(C_0,C_1)$, each giving rise to a {\em hypergeometric
equation} upon multiplying (\ref{dfq2}) by $\zeta(1-\zeta)$:
\be
\zeta(1-\zeta) f'' = [(\alpha+\beta+1) \zeta - \gamma] f' + \alpha\beta f.
\label{hypereq}
\ee
Here $\alpha,\beta,\gamma$ are related to the indices $e_1,e_2,e_3$
at $\zeta=0,1,\infty$ by
\be
\frac1{e_1} = \pm (1-\gamma),\quad
\frac1{e_2} = \pm (\gamma-\alpha-\beta),\quad
\frac1{e_3} = \pm (\alpha-\beta);
\label{hyperabc}
\ee
then $F(\alpha,\beta;\gamma;\zeta)$ and
$(1-\zeta)^\gamma F(\alpha-\gamma+1,\beta-\gamma+1;2-\gamma;\zeta)$
constitute a basis for the solutions of (\ref{dfq2}), where
$F={}_2 F_1$ is the hypergeometric function defined for $|\zeta|<1$ by
\be
F(\alpha,\beta;\gamma;\zeta) := \sum_{n=0}^\infty
 \left[ \prod_{k=0}^{n-1} \frac{(\alpha+k) (\beta+k)} {(\gamma+k)} \right]
 \,\frac{\zeta^n}{n!},
\label{hyperfn}
\ee
and by similar power series in neighborhoods of $\zeta=1$ and
$\zeta=\infty$ (see for instance \cite[9.10 and 9.15]{GR}).
In general, knowing $\sigma$ we may construct and
solve a Schwarzian equation in power series,
albeit series less familiar than ${}_2 F_1$, and numerically
compute the map $X\ra\H/\Gamma$ as the quotient of two solutions.
But once $\Gamma$ is not a triangle group --- that is, when
$X$\/ has more than three elliptic points or positive genus ---
the elliptic points and their orders no longer determine~$\sigma$
but only restrict it to an affine space of finite but positive
dimension.  In general it is a refractory problem to find the
``accessory parameters'' that tell where $\sigma$ lies in that space.
If $\Gamma$ is commensurable with a triangle group $\Gamma'$ then
we obtain $\sigma$ from the quadratic differential on $\H/\Gamma'$
via the correspondence between that curve and $\H/\Gamma$; but this
only applies to Shimura curves associated with the nineteen quaternion
algebras listed by Takeuchi in~\cite{Takeuchi}, including only two
over~$\Q$, the matrix algebra and the algebra ramified at $\{2,3\}$.
One of the advances in the present paper is the computation of $\sigma$
for some arithmetic groups not commensurable with any triangle group.

We now return to the Shimura curves $\Xone$, $\XX(1)$ obtained
from arithmetic groups~$\Gamma=\Gamma(1),\Gamma^*(1)$.
These curves also have a modular interpretation that gives them
the structure of algebraic curves over~$\Q$.  To begin with,
$\Xone$ is the modular curve for principally polarized
abelian surfaces (ppas)~$A$ with an embedding
$O\hookrightarrow{\rm End}(A)$.  (In the classical case
$\O=M_2(\Z)$, corresponding to $\Sigma=\emptyset$, such
an abelian surface is simply the square of an elliptic curve
and we recover the familiar picture of modular curves parametrizing
elliptic ones, but for nonempty $\Sigma$ the surfaces~$A$ are simple
except for those associated to CM points on~$\Xone$; we shall say
more about CM points later.)   The periods of these surfaces
satisfy a linear second-order differential equation which
is a Schwarzian equation for $\H/\Gamma(1)$, usually called
a ``Picard-Fuchs equation'' in this context.
[This generalizes the expression for the periods of elliptic curves
(a.k.a.\ ``complete elliptic integrals'') as ${}_2F_{1}$ values, for
which see e.g.\ \cite[8.113~1.]{GR}.]  The group $\Gamma^*(1)/\Gamma(1)$
acts on $\Xone$ with quotient curve $\XX(1)$.  For each $p\in\Sigma$
there is then an involution $\ww_p\in\Gamma^*(1)/\Gamma(1)$ associated
to the class in~$\Gamma^*(1)/\Gamma(1)$ of elements of~$\O$ of norm~$p$,
and these involutions commute with each other.  (We chose the notation
$\ww_p$ to suggest an analogy with the Atkin-Lehner involutions $w_l$,
which as we shall see have a more direct counterpart in our setting
when $l\notin\Sigma$.)

In terms of abelian surfaces these involutions~$\ww_p$ of~$\Xone$
may be explained as follows.  Let $I_p\subset\O$ consist of the elements
whose norm is divisible by~$p$.  Then $I_p$ is a two-sided prime ideal
of~$\O$, with $\O/I_p\cong\F_{p^2}$ and $I_p^2 = p\O$.  Given
an action of~$\O$ on a ppas~$A$, the kernel of~$I_p$
is a subgroup of~$A$ of size~$p^2$ isotropic under the Weil pairing,
so the quotient surface $A':=A/\ker I_p$ is itself principally
polarized.  Moreover, since $I_p$ is a two-sided ideal, $A'$ inherits
an action of~$\O$.  Thus if $A$ corresponds to some point $P\in\Xone$
then $A'$ corresponds to a point $P'\in\Xone$ determined algebraically
by~$P$\/; that is, we have an algebraic map $\ww_p: P \mapsto P'$
from~$\Xone$ to itself.  Applying this construction to~$A'$ yields
$A/\ker I_p^2 = A/\ker p\O = A / \ker p \cong A$; thus $\ww_p(P')=P$
and $\ww_p$ is indeed an involution.  The quotient curve $\XX(1)$
then parametrizes
% \footnote{
%   Here and later ``parametrizes'' implicitly means ``parametrizes
%   up to isomorphism
%   defined over an algebraic closure''$\!$.
%   }
surfaces~$A$ up to the identification of~$A$ with $A/\ker I$\/ where
$I = \cap_{p\in \Sigma'} I_p = \prod_{p\in \Sigma'} I_p$
for some $\Sigma'\subseteq \Sigma$.

Since $\Xone$, $\XX(1)$ have the structure of algebraic
curves over~$\Q$, they can be regarded as curves over~$\R$.
Now a real structure on any Riemann surface is equivalent
to an anti-holomorphic involution of the surface.  For surfaces
$\H/\Gamma$ uniformized by the upper half-plane, we can give such an
involution by choosing a group $(\Gamma:2)\subset\PGL_2(\R)$
containing $\Gamma$ with index~2 such that
$(\Gamma:2)\not\subset\PSL_2(\R)$.  An element
$({a\;b\atop c\;d})\R^*$ of $\PGL_2(\R)-\PSL_2(\R)$
(i.e.\ with $ad-bc<0$) acts on~$\H$ anti-holomorphically
$z \mapsto (a\bar z + b) / (c\bar z + d)$.  Such a
fractional conjugate-linear transformation
has fixed points on~$\H$ if and only if $a+d=0$, in which case
it is an involution and its fixed points constitute a hyperbolic
line.  Thus $\H/\Gamma$, considered as a curve over~$\R$ using
$\Gamma:2$, has real points if and only if
$(\Gamma:2)-\Gamma$ contains an involution of~$\H$.
The real structures on $\Xone$, $\XX(1)$ are defined by
\bea
(\Gamma(1):2) & := & \O^*/\{\pm1\},
\label{G:2}
\\
(\Gamma^*(1):2) & := &
\{ [a] \in \A^*/\Q^*: a \O = \O a \}.
\label{G*:2}
\eea
That is, compared with (\ref{G},\ref{G*}) we drop the condition
that the norm be positive.  If $\Sigma\neq\emptyset$ then
$\Xone$ has no real points, because if $\Gamma(1):2$
contained an involution $\pm a$ then the characteristic
equation of~$a$ would be $a^2-1=0$ and $\A$ would contain
the zero divisors $a\pm1$.  This is a special case of the
result of~\cite{S3}.  But $\XX(1)$ may have real points.
For instance, we shall see that if $\Sigma=\{2,3\}$ then
$\Gamma^*(1)$ is isomorphic with the triangle group $G_{2,4,6}$.
For general $p,q,r$ with\footnote{
  If $1/p+1/q+1/r$ equals or exceeds~$1$, an analogous situation
  occurs with $\H$ replaced by the complex plane or Riemann sphere.
  }
$1/p+1/q+1/r<1$ we can (and, if $p,q,r$
are distinct, can only) choose $G_{p,q,r}:2$ so that the real
locus of $\H/G_{p,q,r}$ consists of three hyperbolic lines
joining the three elliptic points in pairs, forming a hyperbolic
triangle, with $G_{p,q,r}:2$ generated by hyperbolic reflections
in the triangle's sides; it is this triangle to which the
term ``triangle group'' alludes.

\subsection{The Shimura modular curves $\calX(N)$ and $\XX(N)$
 (with $N$\/~coprime to~$\Sigma$);
the curves $\calXo(N)$ and $\XXo(N)$ and their involution $w_N$}

Now let $l$\/ be a prime not ramified in~$\A$.  Then 
$\A\otimes\Q_l$ and $\O\otimes\Z_l$ are isomorphic with
$M_2(\Q_l)$ and $M_2(\Z_l)$ respectively.  Thus
$(\O\otimes\Q_l)^*_1/\{\pm1\}\cong\PSL_2(\Z_l)$,
with the subscript~$1$ indicating the norm-1 subgroup
as in~(\ref{G}).  We can thus define
congruence subgroups $\Gamma(l)$, $\Gamma_1(l)$, $\Gamma_0(l)$
of~$\Gamma(1)$ just as in the classical case in which
$\Sigma=\emptyset$ and $\Gamma(1)=\PSL_2(\Z)$.
For instance, $\Gamma(l)$ is the normal subgroup
\be
\{ \pm a \in \O^*_+/\{\pm1\} : a \equiv 1\bmod l \}
\label{G(l)}
\ee
of~$\Gamma(1)$, with $\Gamma(1)/\Gamma(l)\cong\PSL_2(\F_l)$;
once we choose an identification of the quotient group
$\Gamma(1)/\Gamma(l)$ with $\PSL_2(\F_l)$ we may define
$\Gamma_0(l)$ as the preimage in $\Gamma(1)$ of the upper
triangular subgroup of $\PSL_2(\F_l)$.
Likewise we have subgroups $\Gamma(l^r)$, $\Gamma_0(l^r)$ etc.,
and even $\Gamma(N)$, $\Gamma_0(N)$ for a positive integer~$N$\/
not divisible by any of the primes of~$\Sigma$.

The quotients of~$\H$ by these subgroups of~$\Gamma(1)$ are then
modular curves covering $\Xone$, which we denote by $\calX(l)$,
$\calXo(l)$, etc.  They parametrize ppas's~$A$ with an
$\O$-action and extra structure: in the case of~$\calX(N)$, a choice
of basis for the $N$-torsion points $A[N]$; in the case of $\calXo(N)$,
a subgroup $G\subset A[N]$ isomorphic with $(\Z/N)^2$ and isotropic
under the Weil pairing.  In the latter case the surface $A'=A/G$
is itself principally polarized and inherits an action of~$\O$
from~$A$, and the image of $A[N]$ in~$A'$ is again a subgroup
$G'\cong(\Z/N)^2$ isotropic under the Weil pairing.  Thus if
we start from some point~$P$\/ on $\calXo(N)$ and associate to
it a pair $(A,G)$ we obtain a new pair $(A',G')$ of the same kind
and a new point $P'\in\calXo(N)$ determined algebraically by~$P$.
Thus we have an algebraic map $w_N: P \mapsto P'$ from~$\calXo(N)$
to itself.  As in the classical case --- in which it is easy to see
that the construction of $A',G'$ from $A,G$\/ amounts to (the square
of) the familiar picture of cyclic subgroups and dual isogenies ---
this $w_N$ is an involution of~$\calXo(N)$ that comes from a trace-zero
element of~$\A$ of norm~$N$\/ whose image in $\A^*/\Q^*$ is an
involution normalizing $\Gamma_0(N)$.

By abuse of terminology we shall say that a pair of points $P,P'$
on~$\Xone$ are ``cyclically $N$-isogenous''\footnote{
  This qualifier ``cyclically'' is needed to exclude cases such
  as the multiplication-by-$m$ map, which as in the case of elliptic
  curves would count as an ``$m^2$-isogeny'' but not a cyclic one.
  }
if they correspond to ppas's $A,A'$ with $A'=A/G$
as above, and call the quotient map $A\ra A/G \cong A'$ a
``cyclic $N$-isogeny''$\!$.  If we regard $P,P'$ as $\Gamma(1)$-orbits
in~$\H$ then they are cyclically $N$-isogenous iff a point in the
first orbit is taken to a point in the second by some $a\in\O$ of
norm~$N$ such that $a\neq m a'$ for any $a'\in\O$ and $m>1$;
since in that case $\bar a$ also satisfies this condition and
acts on~$\H$ as the inverse of~$a$, this relation on $P,P'$ is
symmetric.  Then $\calXo(N)$ parametrizes pairs of $N$-isogenous
points on~$\Xone$, and $w_N$ exchanges the points in such a pair.
The involutions $\ww_p$ on~$\Xone$ lift to the curves $\calX(N)$,
$\calXo(N)$, etc., and commute with $w_N$ on~$\calXo(N)$.

The larger group $\Gamma^*(1)$ likewise has congruence groups
such as $\Gamma^*(N)$, $\Gamma^*_0(N)$, etc., which give rise
to modular curves covering $\XX(1)$ called $\XX(N)$, $\XXo(N)$, etc.
The involution $w_N$ on $\calXo(N)$ descends to an involution on
$\XXo(N)$ which we shall also call $w_N$.  We extend our abuse
of terminology by saying that two points on $\XX(1)$ are ``cyclically
$N$-isogenous'' if they lie under two $N$-isogenous points
of~$\Xone$, and speak of ``$N$-isogenies'' between the equivalence
classes of ppas's parametrized by $\XX(1)$.
One new feature of the congruence subgroups
of~$\Gamma^*(1)$ is that, while $\Gamma^*(N)$ is still normal
in $\Gamma^*(1)$, the quotient group may be larger than
$\PSL_2(\Z/N)$, due to the presence of the $\ww_p$.
For instance if $l\notin\Sigma$ is prime then
$\Gamma^*(1)/\Gamma^*(l)$ is $\PSL_2(\F_l)$ only if all the
primes in~$\Sigma$ are squares modulo~$l$\/; otherwise
the quotient group is $\PGL_2(\F_l)$.  In either case
the index of $\Gamma^*_0(N)$ in $\Gamma^*(1)$, and thus
also the degree of the cover $\XXo(N)/\XX(1)$, is $l+1$.

Since these curves are all defined over~$\Q$, they can again be
regarded as curves over~$\R$ by a suitable choice of $(\Gamma:2)$.
For instance, if $\Gamma = \Gamma(N)$, $\Gamma_1(N)$, $\Gamma_0(N)$
we obtain $(\Gamma:2)$ by adjoining $a\in\O$ of norm $-1$ such
that $a\equiv({1\;\phantom-0\atop0\;-1})\bmod N$\/ under our
identification of $\O/N\O$ with $M_2(\Z/N)$.  Note however that
most of the automorphisms $\PSL_2(\Z/N)$ of $\calX(N)$ do not
commute with $({1\;\phantom-0\atop0\;-1})$ and thus do not
act on $\calX(N)$ regarded as a real curve.  Similar remarks
apply to $\Gamma^*(N)$ etc.

Now fix a prime $l\notin\Sigma$ and consider the sequence of modular
curves $X_r=\calXo(l^r)$ or $X_r=\XXo(l^r)$ ($r=0,1,2,\ldots$).  
The $r$-th curve parametrizes $l^r$-isogenies, which is to say
sequences of $l$-isogenies
\be
A_0 \ra A_1 \ra A_2 \ra \cdots \ra A_n
\label{lisogs}
\ee
such that the composite isogeny $A_{j-1} \ra A_{j+1}$ is a {\em cyclic}
$l^2$-isogeny for each $j$ with $0<j<n$.  Thus for each $m=0,1,\ldots,n$
there are $n+1-m$ maps $\pi_j: X_n \ra X_m$ obtained by
extracting for some $j=0,1,\ldots,n-m$ the cyclic $l^m$-isogeny
$A_j \ra A_{j+m}$ from (\ref{lisogs}).  Each of these maps has degree
$l^{n-m}$, unless $m=0$ when the degree is $(l+1) l^{n-1}$.
In particular we have a tower of maps
\def\pora{\stackrel{\pi_0}{\ra}}
\be
X_n \pora X_{n-1} \pora X_{n-2} \pora \cdots \pora X_2 \pora X_1,
\label{Xotower}
\ee
each map being of degree~$l$.  We observed in~\cite[Prop.~1]{Modtower}
that explicit formulas for $X_1,X_2$, together with their
involutions $w_l,w_{l^2}$ and the map $\pi_0: X_2 \ra X_1$,
suffice to exhibit the entire tower (\ref{Xotower}) explicitly:
For $n\geqs2$ the product map
\be
\pi = \pi_0 \times \pi_1 \times \pi_2 \times \cdots \times \pi_{n-2}:
X_n \ra X_2^{n-1}
\label{prodmap}
\ee
is a 1:1 map from $X_n$ to the set of
$(P_1,P_2,\ldots,P_{n-1}) \in X_2^{n-1}$ such that
\be
\pi_0 \bigl( w_{l^2} (P_j) \bigr) =
w_l \bigl( \pi_0(P_{j+1}) \bigr)
\label{conds}
\ee
for each $j=1,2,\ldots,n-2$.  Here we note that this information
on $X_1,X_2$ is in turn determined by explicit formulas for
$X_0,X_1$, together with the involution $w_l$ and the map
$\pi_0: X_1 \ra X_0$.  Indeed $\pi_1: X_1 \ra X_0$ is then
$\pi_0 \circ w_l$, and the product map
$\pi_0 \times \pi_1: X_2 \ra X_1^2$ identifies $X_2$ with a
curve in $X_1^2$ contained in the locus of
\be
\{ (Q_1,Q_2)\in X_1^2 \,:\, \pi_1(Q_1) = \pi_0(Q_2) \},
\label{locus}
\ee
which decomposes as the union of that curve with the graph
of~$w_l$.\footnote{
  This is where we use the hypothesis that $l$\/ is prime.  The
  description of $X_n$ in (\ref{prodmap},\ref{conds}) holds even
  for composite~$l$, but the description of $X_2$ in terms of
  $X_1$ does not, because then (\ref{locus}) has other components.
  }
This determines $X_2$ and the projections $\pi_j: X_2 \ra X_1$
($j=0,1$); the involution $w_{l^2}$ is
\be
(Q_1,Q_2) \lra (w_l Q_2, w_l Q_1).
\label{w(l^2)}
\ee
Thus the equations we shall exhibit for certain choices of~$\A$
and~$l$\/ suffice to determine explicit formulas for towers
of Shimura modular curves $\calXo(l^r)$, $\XXo(l^r)$, towers whose
reduction at any prime $l'\notin\Sigma\cup\{l\}$ is known to be
asymptotically optimal over the field of ${l'}^2$ elements
\cite{I3,TVZ}.

\subsection{Complex-multiplication (CM) and supersingular points
on Shimura curves}

Let $F$\/ be a quadratic imaginary field, and let $O_F$ be its
ring of integers.  Assume that none of the primes of~$\Sigma$
split in~$F$.  Then $F$\/ embeds in~$\A$ (in many ways), and
$O_F$ embeds in~$\O$.  For any embedding $\iota:F\hookrightarrow\A$,
the image of $F^*$ in $\A^*/\Q^*$ then has a unique fixed point
on $\H$; the orbit of this point under $\Gamma(1)$, or under any
other congruence subgroup $\Gamma\subset\A^*/\Q^*$, is then a
{\em CM point}\/ on the Shimura curve $\H/\Gamma$.  In particular,
on $\Xone$ such a point parametrizes a ppas with extra endomorphisms
by~$\iota(F)\cap\O$.  For instance if $\iota(F)\cap\O=\iota(O_F)$ then this
ppas is a product of elliptic curves each with complex multiplication
by $O_F$ (but not in the product polarization).  In general
$\iota^{-1}(\iota(F)\cap\O)$ is called the {\em CM ring} of the
CM~point on~$\Xone$.  Embeddings conjugate by $\Gamma(1)$ yield the
same point on $\Xone$, and for each order $O\subset F$\/ there are
finitely many embeddings up to conjugacy, and thus finitely many
CM~points on~$\Xone$ with CM~ring~$O$\/; in fact their number is
just the class number of~$O$.  In~\cite{S2} Shimura already
showed that all points with the same CM~ring are Galois conjugate
over~$\Q$, from which it follows that a CM~point is rational if
and only if its CM~ring has unique factorization.

Thus far the description is completely analogous to the theory
of complex multiplication for $j$-invariants of elliptic curves.
But when $\Sigma\neq\emptyset$ a new phenomenon arises: CM~points
on the quotient curve $\XX(1)$ may be rational even when their
preimages on~$\Xone$ are not.  For instance, a point with CM~ring
$O_F$ is rational on~$\XX(1)$ if and only if the class group of~$F$\/
is generated by the classes of ideals $I\subset O_F$\/ such that
$I^2$ is the principal ideal $(p)$ for some rational prime $p\in\Sigma$.
This has the amusing consequence that when $\Sigma=\{2,3\}$ the number
of rational CM~points on $\XX(1)$ is more than twice the number of
rational CM~points on the classical modular curve $X(1)$.
[Curiously, already in the classical setting $X(1)$
does not hold the record: it has 13 rational CM points,
whilst $X_0^*(6)=X_0(6)/\langle w_2,w_3 \rangle$ has~14.
The reason again is fields~$F$\/ with nontrivial class group
generated by square roots of the ideals $(2)$ or $(3)$, though with
a few small exceptions both 2 and~3 must ramify in~$F$.  In the
$\XX(1)$ setting the primes of~$\Sigma$ are allowed to be inert
as well, which makes the list considerably longer.]  In fact for
each of the first four cases $\Sigma = \{2,3\}, \{2,5\}, \{2,7\},
\{3,5\}$ we find more rational CM points than on any classical
modular curve.

A major aim of this paper is computation of the coordinates of
these points.  We must first list all possible~$O$.  The class
number of~$O$, and thus of~$F$, must be a power of~2 no greater
than $2^{\#\Sigma}$.  In each of our cases, $\#\Sigma=2$,
so $F$\/ has class number at most~4 and we may refer to the
list of imaginary quadratic number fields with class group
$(\Z/2)^r$ ($r=0,1,2$), proved complete by Arno~\cite{Arno}.\footnote{
  It might be possible to avoid that difficult proof for our
  application, since we are only concerned with fields whose class
  group is accounted for by ramified primes in a given set~$\Sigma$,
  and it may be possible to provably list them all using the arithmetic
  of CM points on either classical or Shimura modular curves, as in
  Heegner's proof that $\Q(\sqrt{-163}\,)$ is the last quadratic
  imaginary field of class number~1.
  }
Given $F$\/ we easily find all possible~$O$, and imbed each into~$\O$
by finding $a\in\O$ such that $(a-\bar a)^2 = \mathop{\rm disc}(O)$.
This gives us the CM point on~$\H$.  But we want its coordinates
on the Shimura curve $\H/\Gamma^*(1)$ as rational numbers.  Actually
only one coordinate is needed because $\XX(1)$ has genus~0 for
each of our four $\Sigma$.  We recover the coordinate as a real
number using our Schwarzian uniformization of $\XX(1)$ by~$\H$.
(Of course a coordinate on $\Pr^1$ is only defined up to $\PGL_2(\Q)$,
but in each case we choose a coordinate once and for all by specifying
it on the CM~points.)  We then recognize that number as a rational
number from its continued fraction expansion, and verify that the
putative rational coordinate not only agrees with our computations
to as many digits as we want but also satisfies various arithmetic
conditions such as those described later in this section.

Of course this is not fully satisfactory; we do not know how to
prove that, for instance,
$t = 13^2 67^2 109^2 139^2 157^2 163 / 2^{10} 5^6 11^6 17^6$
(see Tables 1,2 below) is the CM point of discriminant $-163$
on the curve $\XX(1)$ associated with the algebra ramified
at $\{2,3\}$.  But we can prove that above half of our numbers
are correct, again using the modular curves $\XX(l)$ and their
involutions $w_l$ for small~$l$.  This is because CM points behave
well under isogenies: any point isogenous to a CM point is itself CM,
and moreover a point on $\Xone$ or $\XX(1)$ is CM if and only if it
admits a cyclic $d$-isogeny to itself for some $d>1$.  Once we have
formulas for $\XXo(l)$ and $w_l$ we may compute all points cyclically
$l$-isogenous either with an already known CM points or with themselves.
The discriminant of a new rational CM point can then be determined
either by arithmetic tests or by identifying it with a real CM~point
to low precision.

The classical theory of {\em supersingular points} also largely
carries over to the Shimura setting.  We may use the fact that
the ppas parametrized by a CM point has extra endomorphisms to
define CM~points of Shimura curves algebraically, and thus in
any characteristic $\notin\Sigma$.  In positive characteristic
$p\notin\Sigma$, any CM point is defined over some finite field, and
conversely every $\overline{\F_{\!p}}$-point of a Shimura curve is CM.
All but finitely many of these parametrize ppas's whose endomorphism
ring has $\Z$-rank~8; the exceptional points, all defined over
$\F_{\!p^2}$, yield rank~16, and are called {\em supersingular},
all other $\overline{\F_{\!p}}$-points being {\em ordinary}.  One may
choose coordinates on~$\Xone$ (or $\XX(1)$) such that a CM~point
in characteristic~zero reduces mod~$p$ to a ordinary point 
if $p$ splits in the CM~field, and to a supersingular point
otherwise.  Conversely each ordinary point mod~$p$ lifts to a unique
CM~point (cf.~\cite{Deuring} for the classical case).  This means
that if two CM~points with different CM~fields have the same reduction
mod~$p$, their common reduction is supersingular, and then as
in~\cite{GZ} there is an upper bound on~$p$ proportional to the
product of the two CM~discriminants.  So for instance if
$\XX(1)\cong\Pr^1$ then the difference between the coordinates
of two rational CM~points is a product of small primes.  This
remains the case, for similar reasons, even for distinct CM~points
with the same CM~field, and may be checked from the tables
of rational CM~points in this paper.  The preimages of the
supersingular points on modular covers such as $\calXo(l)$ yield
enough $\F_{\!p^2}$-rational points on these curves to attain
the Drinfeld-Vl\u{a}du\c{t} bound~\cite{I3}; these curves are
thus ``asymptotically optimal'' over $\F_{\!p^2}$.  Asymptotically
optimal curves over $\F_{\!p^{2f}}$ ($f>1$) likewise come from
Shimura curves associated to quaternion algebras over totally
real number fields with a prime of residue field~$\F_{\!p^f}$.

In the case of residue field $\F_{\!p}$ (so in particular for
quaternion algebras over~$\Q$) Ihara~\cite{I2} found a
remarkable connection between the hyperbolic uniformization
of a Shimura curve~$\calX=\H/\Gamma$ and the supersingular points of
its reduction mod~$p$.  We give his result in the case that
$\calX$ has genus~0, because we will only apply it to such
curves and the result can be stated in an equivalent and elementary
form (though the {\em proof}\/ is still far from elementary).
Since we are working over $\F_{\!p}$, we may identify any curve
of genus~0 with $\Pr^1$, and choose a coordinate (degree-1 function)
$t$ on $\Pr^1$ such that $t=\infty$ is an elliptic point.  Let
$t_i$ be the coordinates of the remaining elliptic points.

First, the hyperbolic area of the curve controls the number of
points, which is approximately $\frac12(p+1){\rm Area}(\calX)$ ---
``approximately'' because $\frac12(p+1){\rm Area}(\calX)$ is not
the number of points but their total mass.  The mass of a non-elliptic
supersingular point is~$1$, but an elliptic point with stabilizer~$G$\/
has mass $1/\#G$.  If the elliptic point mod~$p$ is the reduction of
only one elliptic point on $\H/\Gamma$ (which, for curves
coming from quaternion algebras over~$\Q$, is always the
case once $p>3$), then its stabilizer is $\Z/e\Z$ and its
mass is $1/e$ where $e$ is the index of that elliptic point.
[The mass formula also holds for $\calX$ of arbitrary genus,
and for general residue fields provided $p$ is replaced by
the size of the field.]  Let $d$\/ be the number of non-elliptic
supersingular points, and choose a Schwarzian equation (\ref{dfq2})
with at most regular singularities at $t=\infty,t_i$ and no other
singularities.  Then the supersingular points are determined uniquely
by the condition that their $t$-coordinates are the roots of a
polynomial $P(t)$ of degree~$d$\/ such that for some $r_i\in\Q$
the algebraic function $\prod_i (t-t_i)^{r_i} \cdot P(t)$
is a solution of the Schwarzian differential equation~(\ref{dfq2})!
For instance~\cite[4.3]{I2}, if $\Gamma$ is a triangle group
we may choose $t_i=0,1$, and then $P(t)$ is a finite hypergeometric
series mod~$p$.

Given $t_0\in\Q$ we may then test whether $t_0$ is ordinary or
supersingular mod~$p$ for each small~$p$.  If $t_0$ is a CM~point
with CM~field then its reduction is ordinary if $p$ splits in~$F$,
supersingular otherwise.  When we have obtained $t_0$ as a good
rational approximation to a rational CM point, but could not prove
it correct, we checked for many $p$ whether $t_0$ is ordinary or
supersingular mod~$p$; when each prime behaves as expected from its
behavior in~$F$, we say that $t_0$ has ``passed the supersingular
test'' modulo those primes~$p$.

\section{The case $\Sigma = \{2,3\}$}
\subsection{The quaternion algebra and the curves $\Xone$, $\XX(1)$}

For this section we let $\A$ be the quaternion algebra ramified
at $\{2,3\}$.  This algebra is generated over~$\Q$ by elements
$b,c$ satisfying
\be
b^2=2,\ \ c^2=-3,\ \ bc=-cb.
\label{bcrel}
\ee
The conjugation of~$\A$ fixes~1 and takes $b,c,bc$ to $-b,-c,-bc$;
thus for any element $\alpha = \alpha_1 + \alpha_2 b + \alpha_3 c
+ \alpha_4 bc \in \A$ the conjugate and norm of $\alpha$ are given by
\be
\bar\alpha = \alpha_1 - \alpha_2 b - \alpha_3 c - \alpha_4 bc,
\quad
\Nm(\alpha) = \alpha_1^2 - 2\alpha_2^2 + 3\alpha_3^2 - 6\alpha_4^2.
\ee
Since $\A$ is indefinite, all its maximal orders are conjugate;
let $\O$ be the maximal order generated by~$b$ and $(1+c)/2$.  Then
$\Gamma^*(1)$ contains $\Gamma(1)$ with index $2^{\#\Sigma}=4$,
and consists of the classes mod~$\Q^*$ of elements of~$\O$ of norm
1, 2, 3, or~6.  In row~II of Table~3 of~\cite{Takeuchi} (p.208)
we find that $\Gamma^*(1)$ is isomorphic with the triangle group
\be
G_{2,4,6} :=
\langle s_2,s_4,s_6 | s_2^2 = s_4^4 = s_6^6 = s_2 s_4 s_6 = 1 \rangle
.
\label{G246}
\ee
Indeed we find that $\Gamma^*(1)$ contains elements
\be
s_2 = [bc+2c],\quad s_4 = [(2+b)(1+c)],\quad s_6 = [3+c]
\label{s246}
\ee
[NB $(2+b)(1+c), 3+c \in 2\O$] of orders $2,4,6$ with $s_2 s_4 s_6 = 1$.
The subgroup of~$\Gamma^*(1)$ generated by these elements is thus
isomorphic with $G_{2,4,6}$.  But a hyperbolic triangle group
cannot be isomorphic with a proper subgroup (since the areas
of the quotients of $\H$ by the group and its subgroup are equal),
so $\Gamma^*(1)$ is generated by $s_2, s_4, s_6$.  Note that these
generators have norms $6,2,3$ mod $(\Q^*)^2$, and thus represent
the three nontrivial cosets of~$\Gamma^*(1)$ in $\O^*/\{\pm1\}$.

Since $\Gamma^*(1)$ is a triangle group, $\XX(1)$ is a curve of
genus~0.  Moreover $\XX(1)$ has $\Q$-rational points (e.g.\ the three
elliptic points, each of which must be rational because it is
the only one of its index), so $\XX(1)\cong\Pr^1$ over~$\Q$.
Let $t$\/ be a rational coordinate on that curve (i.e.\
a rational function of degree~1).  In general a rational coordinate
on $\Pr^1$ is determined only up to the $\PGL_2$ action on~$\Pr^1$,
but can be specified uniquely by prescribing its values at three points.
In our case $\XX(1)$ has three distinguished points, namely the
elliptic points of orders $2,4,6$; we fix $t$\/ by requiring that
it assume the values $0,1,\infty$ respectively at those three points.

None of $s_2,s_4,s_6$ is contained in $\Gamma(1)$.  Hence the
$(\Z/2)^2$ cover $\Xone/\XX(1)$ is ramified at all three elliptic
points.  Thus $s_2$ lies under two points of~$\Xone$ with trivial
stabilizer, while $s_4$ lies under two points of index~2 and
$s_6$ under two points of index~3.  By either the Riemann-Hurwitz
formula or from~(\ref{genus}) we see that $\Xone$ has genus~0.
This and the orders $2,2,3,3$ of the elliptic points do not completely
specify $\Gamma(1)$ up to conjugacy in $\PSL_2(\R)$: to do that we
also need the cross-ratio of the four elliptic points.  Fortunately
this cross-ratio is determined by the existence of the cover
$\Xone\ra\XX(1)$, or equivalently of an involution $s_4$ on~$\Xone$
that fixes the two order-2 points and switches the order-3 points.
This forces the pairs of order-2 and order-3 points to have a
cross-ratio of~$-1$, or to ``divide each other harmonically''
as the Greek geometers would say.  The function field of~$\Xone$
is generated by the square roots of $c_0 t$ and~$c_1(t-1)$ for some
$c_0,c_1\in\Q^*/{\Q^*}^2$, but we do not yet know which multipliers
$c_0,c_1$ are appropriate.  If both $c_0,c_1$ were~1 then $\Xone$
would be a rational curve with coordinate~$u$ with
$t=((u^2+1)/2u)^2 = 1 + ((u^2-1)/2u)^2$,
the familiar parametrization of Pythagorean triples.  The elliptic
points of order~2 and~3 would then be at $u=\pm1$ and $u=0,\infty$.
However it will turn out that the correct choices are $c_0=-1, c_1=3$,
and thus that $\Xone$ is the conic with equation
\be
X^2+Y^2+3Z^2=0
\label{X1.6}
\ee
and no rational points even over~$\R$.  [That $\Xone$ is the conic
(\ref{X1.6}) is announced in \cite[p.279]{Ku} and attributed to Ihara;
that there are no real points on the Shimura curve $\Xone$ associated
to any indefinite quaternion algebra over~$\Q$ other than $M_2(\Q)$
was already shown by Shimura~\cite{S3}.  The equation (\ref{X1.6})
for~$\Xone$ does not uniquely determine $c_0,c_1$, but the local
methods of~\cite{Ku} could probably supply that information as well.]

\subsection{Shimura modular curves $\XXo(l)$ and $\calX(l)$
for $l=5,7,13$}

Let $l$ be a prime other than the primes $2,3$ of~$\Sigma$.
We determine the genus of the curve $\XXo(l)$ using the formula
(\ref{genus}).  Being a cover of $\XX(1)$ of degree $l+1$, the
curve $\XXo(l)$ has normalized hyperbolic area $(l+1)/12$.
It has $1+(-6/l)$ elliptic points of order~2, $1+(-1/l)$
elliptic points of order~4, and $1+(-3/l)$ elliptic points
of order~6.  This is a consequence of our computation of
$s_2,s_4,s_6$, which lift to elements of~$\A$ that generate
subfields isomorphic with $\Q(\sqrt{-6}\,)$, $\Q(\sqrt{-1}\,)$,
and $\Q(\sqrt{-3}\,)$.  Actually the orders $2,4,6$ of the
elliptic points suffice.  Consider the images of $s_2,s_4,s_6$
  in the Galois group ($\subseteq\PGL_2(\F_l)$) of the cover
  $\XXo(l)/\XX(1)$, and the cycle
structures of their actions on the $l+1$ points of $\Pr^1(\F_l)$.
These images $\sigma_2,\sigma_4,\sigma_6$ are group elements of order
$2,4,6$.  For 4 and~6, the order determines the conjugacy class, which
joins as many of the points of $\Pr^1(\F_l)$ as possible in cycles
of length~4 or~6 respectively and leaves any remaining points fixed;
the number of fixed points is two or none according to the
residue of~$l$ mod~4 or~6.  For $\sigma_2$ there are two conjugacy
classes in $\PGL_2(\F_l)$, one with two fixed points and
the other with none, but the choice is determined by the condition
that the genus $g(\XXo(l))$ be an integer, or equivalently by the
requirement that the signs of $\sigma_2,\sigma_4,\sigma_6$ considered
as permutations of $\Pr^1(\F_l)$ be consistent with $s_2 s_4 s_6=1$.
We readily check that this means that the image of $s_2$ has two
fixed points if and only if $(-6/l)=+1$, as claimed.  From
(\ref{genus}) we conclude that
\be
g(\XXo(l)) = \frac1{24} \left[ l
- 6 \left(\frac{-6}{l}\right)
- 9 \left(\frac{-1}{l}\right)
- 10 \left(\frac{-3}{l}\right)
\right] .
\label{genus(6)}
\ee
We tabulate this for $l<50$:

\vspace*{2ex}
 
\setlength{\tabcolsep}{4pt}
\centerline{
\begin{tabular}{c|ccccccccccccc}
$l$ &
5 & 7 & 11 & 13 & 17 & 19 & 23 & 29 & 31 & 37 & 41 & 43 & 47 \\ \hline
$g(\XXo(l))$ & 
0 & 0 &  1 &  0 &  1 &  1 &  2 &  1 &  1 &  1 &  2 &  2 &  3
\end{tabular}
}
\setlength{\tabcolsep}{1.4pt}

It so happens that in the first seven cases $g(\XXo(l))$ coincides
with the genus of the classical modular curve X$_0(l)$, but of course
this cannot go on forever because the latter genus is $l/12 + O(1)$
while the former is only $l/24 + O(1)$, and indeed $g(\XXo(l))$ is
smaller for all $l>23$.  Still, as with X$_0(l)$, we find that
$\XXo(l)$ has genus~0 for $l=5,7,13$, but not for $l=11$ or any $l>13$.
For the three genus-0 cases we shall use the ramification behavior
of the cover $\XXo(l)/\XX(1)$ to find an explicit rational function
of degree $l+1$ on~$\Pr^1$ that realizes that cover and determine the
involution $w_l$.

Now for any $l>3$\/ the solution of $\sigma_2 \sigma_4 \sigma_6 = 1$
in elements $\sigma_2,\sigma_4,\sigma_6$ of orders $2,4,6$ in
$\PGL_2(\F_l)$ is unique up to conjugation in that group.
Thus we know from the general theory of~\cite{Mat} that the
cover $\XXo(l)/\XX(1)$ is determined by its Galois group and
ramification data.  Unfortunately the proof of this fact does
not readily yield an efficient computation of the cover;
for instance the Riemann existence theorem for Riemann surfaces
is an essential ingredient.  We use a method for finding the rational
function $t: \XXo(l) \ra \XX(1)$ explicitly that amounts to
solving for its coefficients, using the cycle structures
of $\sigma_2,\sigma_4,\sigma_6$ to obtain algebraic conditions.
In effect these conditions are the shape of the divisors
$(t)_0$, $(t)_1$, $(t)_\infty$.  But a rational function
satisfying these conditions is not in general known to have
the right Galois group: all we know is that the monodromy
elements around $0,1,\infty$ have the right cycle structures
in the symmetric group S$_{l+1}$.  Thus we obtain several
candidate functions, only one of which has Galois group
$\PGL_2(\F_l)$ (or $\PSL_2(\F_l)$ if $l\equiv1\bmod24$).
Fortunately for $l=5,7$ we can exclude the impostors
by inspection, and for $l=13$ the computation has already
been done for us.

\underline{l=5}.
Here the cycle structures of $s_2,s_4,s_6$ are 2211,~411,~6.
Curiously if the identity in the symmetric group~$S_6$ is written
as the product of three permutations $\sigma_2,\sigma_4,\sigma_6$
with these cycle structures then they can never generate all
of~$S_6$.  This can be seen by considering their images
$\sigma'_2,\sigma'_4,\sigma'_6$ under an outer automorphism of~$S_6$:
these have cycle structures 2211,~411,~321, and thus have too many
cycles to generate a transitive subgroup (if two permutations
of~$n$ letters generate a transitive subgroup of~$S_n$ then they
and their product together have at most $n+2$ cycles).  It turns
out that the subgroup generated by $\sigma'_2,\sigma'_4,\sigma'_6$
can be either $A_4\times S_2$ or the point stabilizer $S_5$.
In the former case $\sigma_2,\sigma_4,\sigma_6$ generate a transitive
but imprimitive subgroup of~$S_6$: the six letters are partitioned
into three pairs, and the group consists of all permutations that
respect this partition and permute the pairs cyclically.  In the
latter case $\sigma_2,\sigma_4,\sigma_6$ generate $\PGL_2(\F_5)$;
this is the case we are interested in.  In each of the two cases the
triple $(\sigma_2,\sigma_4,\sigma_6)$ is determined uniquely up to
conjugation in the subgroup of~$S_6$ generated by the $\sigma$'s, each
of which is in a rational conjugacy class in the sense of~\cite{Mat}.
Thus each case corresponds to a unique degree-6 cover $\Pr^1\ra\Pr^1$
defined over~$\Q$.  We shall determine both covers.

Let $t$ be a rational function on~$\Pr^1$ ramified only above
$t=0,1,\infty$ with cycle structures 2211,~411,~6.  Choose
a rational coordinate~$x$ on~$\Pr^1$ such that $x=\infty$ is
the sextuple pole of~$t$\/ and $x=0$ is the quadruple zero
of $t-1$; this determines $x$ up to scaling.  Then $t$
is a polynomial of degree~6 in~$x$ with two double roots such that
$t\equiv1\bmod x^4$.  The double roots are necessarily the roots
of the quadratic polynomial $x^{-3} dt/dx$.  Thus $t$ is a
polynomial of the form $c_6 x^6 + c_5 x^5 + c_4 x^4 + 1$
divisible by $6 c_6 x^2 + 5 c_5 x + 4 c_4$.  We readily compute
that there are two possibilities for $c_4,c_5,c_6$ up to scaling
$(c_4,c_5,c_6)\ra(\lambda^4 c_4,\lambda^5 c_5,\lambda^6 c_6)$.
One possibility gives $t = 2x^6 - 3x^4 + 1 = (x^2-1)^2 (2x^2+1)$;
being symmetric under $x\lra -x$ this must be the imprimitive solution.
Thus the remaining possibility must give the $\PGL_2(\F_5)$ cover
$\XXo(5) / \XX(1)$.  The following choice
of scaling of $x=x_5$ seems simplest:
$$
t = 540 x^6 + 324 x^5 + 135 x^4 + 1
$$\vspace*{-4.5ex}\be\label{X5}\ee\vspace*{-4.5ex}$$
= 1 + 27x^4(20x^2+12x+5) = (15x^2-6x+1)(6x^2+3x+1)^2 .
$$
The elliptic points of order~2 and~4 on $\XXo(5)$ are the simple
zeros of $t$ and $t-1$ respectively, i.e.\ the roots of
$15x^2-6x+1$ and $20x^2+12x+5$.  The involution $w_5$ switches
each elliptic point with the other elliptic point of
the same order; this suffices to determine~$w_5$.  The fact
that two pairs of points on~$\Pr^1$ switched by an involution
of~$\Pr^1$ determine the involution is well-known, but we have
not found in the literature an explicit formula for doing this.
Since we shall need this result on several occasion we give it
in an Appendix as Proposition~A.  Using that formula~(\ref{propA}),
we find that
\be
w_5(x) = \frac{42 - 55x} {55 + 300x}.
\label{w5}
\ee

\underline{l=7}.
This time $s_2,s_4,s_6$ have cycle structures 22211,~44,~611.
Again there are several ways to get the identity permutation on
8 letters as a product of three permutations with these
cycle structures, none of which generate the full symmetric
group $S_8$.  There are two ways to get the imprimitive group
$2^4:S_4$; the corresponding covers are obtained from the $S_4$
cover $t=4\xi^3-3\xi^4$ by taking $\xi=x^2+\xi_0$ where $\xi_0$
is either root of the quadratic $3\xi^2+2\xi+1 = (1-t)/(\xi-1)^2$.
The remaining solution corresponds to our $\PGL_2(\F_7)$ cover.

To find that cover, let $t$ be a rational function on~$\Pr^1$
ramified only above $t=0,1,\infty$ with cycle structures 2211,~411,~6,
and choose a rational coordinate~$x$ on~$\Pr^1$ such that $x=\infty$ is
the sextuple pole of~$t$.  This determines $x$ up to an affine linear
transformation.  Then there is a cubic polynomial~$P$\/ and quadratic
relatively prime polynomials $Q_1,Q_2,Q_3$ in~$x$ such that
$t=P^2 Q_1/Q_3 = 1 + Q_2^4/Q_3$, i.e.\ such that
$P^2 Q_1 - Q_2^4$ is quadratic.  Equivalently, the Taylor expansion
of $Q_2^2 / \sqrt{Q_1}$ about $x=\infty$ should have vanishing $x^{-1}$
and $x^{-2}$ coefficients, and then $R(x)$ is obtained by truncating
that Taylor expansion after its constant term.  We assume without
loss of generality that $Q_1,Q_2$ are monic.  By translating $x$
(a.k.a.\ ``completing the square'') we may assume that $Q_1$ is
of the form $x^2+\alpha$.  If the same were true of $Q_2$ then $t$\/
would be a rational function of~$x^2$ and we would have an imprimitive
cover.  Thus the constant coefficient of $Q_2$ is nonzero, and by
scaling $x$ we may take $Q_2 = x^2+x+\beta$.  We then set the
$x^{-1},x^{-2}$ coefficients of of $Q_2^2 / \sqrt{Q_1}$ to zero,
obtaining the equations
\be
3\alpha^2-8\alpha\beta+8\beta^2-4\alpha = 3\alpha^2-4\alpha\beta = 0.
\label{eqs7}
\ee
Thus either $\alpha=0$ or $\alpha=4\beta/3$.  The first option
yields $\beta=0$ which fails because then $Q_1,Q_2$ have the
common factor~$x$.  The second option yields $\beta=0$, which
again fails for the same reason, but also $\beta=2$ which succeeds.
Substituting $-(2x+1)/3$ for $x$ to reduce the coefficients we then
find:
$$
t = - \frac
  {(4x^2 + 4x + 25) (2x^3 - 3x^2 + 12x - 2)^2}
  {108(7x^2 - 8x + 37)}
$$\vspace*{-3ex}\be\label{X7(6)}\ee\vspace*{-3ex}$$
  = 1 - \frac{(2x^2 - x + 8)^4} {108(7x^2 - 8x + 37)}.
$$
The elliptic points of order~2 and~6 on $\XXo(7)$ are respectively
the simple zeros and poles of~$t$, i.e.\ the roots of
$4x^2+4x+25$ and $7x^2-38x+7$.  The involution $w_7$ is again
by the fact that it switches each elliptic point with the other
elliptic point of the same order: it is 
\be
w_7(x) = \frac{116 - 9x} {9 + 20x}.
\label{w7(6)}
\ee

\underline{l=13}.
Here the cycle structures are $2^7$, 44411, 6611.  The computation
of the degree-14 map is of course much more complicated than for
the maps of degrees $6,8$ for $l=5,7$.  Fortunately this computation
was already done in \cite[\S4]{MM} (a paper concerned not with Shimura
modular curves but with examples of rigid $\PSL_2(\F_p)$ covers
of the line).  There we find that there is a coordinate $x=x_{13}$
on~$\XXo(13)$ for which
$$
t = 1 - \frac{27}{4} \frac
   {(x^2+36) (x^3+x^2+35x+27)^4}
   {(7x^2+2x+247) (x^2+39)^6}
$$\vspace*{-4.5ex}\be\label{X13}\ee\vspace*{-1.5ex}$$
  = \frac
   {(x^7-50x^6+63x^5-5040x^4+783x^3-168426x^2-6831x-1864404)^2}
   {4 (7x^2+2x+247) (x^2+39)^6}.
$$
The elliptic points of order~4 and~6 on $\XXo(13)$ are respectively
the simple zeros and poles of~$t-1$, i.e.\ the roots of $x^2+36$
and $7x^2+2x+247$.  Once more we use (\ref{propA}) to find the
involution from the fact that it switches each elliptic point
with the other elliptic point of the same order:
\be
w_{13}(x) = \frac{5x+72} {2x-5}.
\label{w13}
\ee

{}From an equation for $\XX(l)$ and the rational map~$t$ on that curve
we recover $\calXo(l)$ by adjoining square roots of $c_0 t$ and
$c_1 (t-1)$.  For each of our three cases $l=5,7,13$ the resulting
curve has genus~1, and its Jacobian is an elliptic curve of
conductor~$6l$ --- but only if we choose $c_0,c_1$ that give the
correct quadratic twist.  For $l=5$, $l=7$, $l=13$ it turns out that
we must take a square root of $3t(1-t)$, $-t$, $3(t-1)$ respectively.
Fortunately these are consistent and we obtain \hbox{$c_0=-1$} and
\hbox{$c_1=3$} as promised.  The resulting curves
$\calXo(5),\calXo(7),\calXo(13)$
have no rational or even real points (because this is already true
of the curve $\Xone$ which they all cover); their Jacobians are
the curves numbered 30F,\,42C,\,78B in the Antwerp tables in~\cite{BK}
compiled by Tingley et al., and and 30-A6,\,\hbox{42-A3,}\,\hbox{78-A2}
in Cremona~\cite{Cr}.

\subsection{Supersingular points on~$\XX(1)\bmod l$}

We have noted that Ihara's description of supersingular points on
Shimura curves is particularly simple in the case of a triangle group:
the non-elliptic supersingular points are roots of a hypergeometric
polynomial, and the elliptic points are CM in characteristic zero
so the Deuring test determines whether each one is supersingular or not.

In our case, The elliptic points $t=0$, $t=1$, $t=\infty$
are supersingular mod~$l$\/ if and only iff $l$\/ is inert
in $\Q(\sqrt{-6}\,)$, $\Q(\sqrt{-1}\,)$, $\Q(\sqrt{-3}\,)$
respectively, i.e.\ iff $-6$, $-1$, $-3$ is a quadratic
nonresidue of~$l$.  Thus the status of all three elliptic
points depends on $l\bmod 24$, as shown in the next table:

\vspace*{2ex}

\setlength{\tabcolsep}{4pt}
\centerline{
\def\bt{$\bullet$}
\def\hi#1{\raise1.5ex\hbox{$#1$}}
\begin{tabular}{c|c||c|c|c|c|c|c|c|c}
        &     & \multicolumn{8}{c}{$l$ mod 24}     \\
\hi t   &\hi e& 1 & 5 & 7 & 11 & 13 & 17 & 19 & 23 \\ \hline
    0   &  2  &   &   &   &    &\bt &\bt &\bt &\bt \\
    1   &  4  &   &   &\bt&\bt &    &    &\bt &\bt \\
$\infty$&  6  &   &\bt&   &\bt &    &\bt &    &\bt \\
\end{tabular}
}
\setlength{\tabcolsep}{1.4pt}

\vspace*{2ex}

\noindent
(bullets mark elliptic points with supersingular reduction).
This could also be obtained from the total mass $(l+1)/24$
of supersingular points, together with the fact that
the contribution to this mass of the non-elliptic points
is integral: in each column the table shows the unique subset of
$1/2,1/4,1/6$ whose sum is congruent to $(l+1)/24$ mod~1.
The hypergeometric polynomial whose roots are the non-elliptic
supersingular points has degree $\lfloor l/24 \rfloor$, and
depends on $l\bmod24$ as follows:
\be
\def\OO{\phantom{|\atop|}}
\cases{
F(\frac{1}{24},\frac{5}{24};\frac12;t),\OO& if $l\equiv1$ or~5 mod~24;
 \cr
F(\frac{7}{24},\frac{11}{24};\frac12;t),\OO& if $l\equiv7$ or~11 mod~24;
 \cr
F(\frac{13}{24},\frac{17}{24};\frac32;t),\OO& if $l\equiv13$ or~17 mod~24;
 \cr
F(\frac{19}{24},\frac{23}{24};\frac32;t),\OO& if $l\equiv19$ or~23 mod~24.\cr
}
\label{ssing}
\ee

For example, for $l=163 (\equiv19\bmod24)$ we find
\bea
F(\frac{19}{24},\frac{23}{24};\frac32;t) &=&
43t^6 + 89t^5 + 97t^4 + 52t^3 + 149t^2 + 132t + 1 \nonumber \\
&& =
(t+76) (t+78) (t+92) (t+127) (t^2+65t+74)
\label{163}
\eea
in characteristic 163, so the supersingular points mod~163 are
$0,1$, and the roots of~(\ref{163}) in $\F_{163^2}$.

\subsection{CM points on~$\XX(1)$ via $\XXo(l)$ and $w_l$}

We noted already that the elliptic points $t=0,1,\infty$ on~$\XX(1)$
are CM points, with discriminants $-3,-4,-24$.  Using our formulas
for $\calXo(l)$ and $w_l$ ($l=5,7,13$) we can obtain fourteen further
CM points: three points isogenous to one of the elliptic CM points, and
eleven more points cyclically isogenous to themselves.  This accounts
for all but ten of the 27 rational CM points on~$\XX(1)$.  The
discriminants of the three new points isogenous to $t=1$ or $t=\infty$
are determined by the isogenies' degrees.  The discriminants of the
self-isogenous points can be surmised by testing them for supersingular
reduction at small primes: in each case only one discriminant small
enough to admit a self-isogeny of that degree has the correct quadratic
character at the first few primes, which is then confirmed by extending
the test to all primes up to~200.

On $\XXo(5)$ the image of $x_5=\infty$ under $w_5$ is $-11/60$, which
yields the CM point $t=152881/138240$; likewise from $w_5(0)=42/55$
we recover the point $421850521/1771561$.  These CM points are
5-isogenous with the elliptic points \hbox{$t=\infty$}, $t=1$
respectively, and thus have discriminants $-3\cdot 5^2$ and
$-4\cdot 5^2$.  Similarly on $\XXo(7)$ we have $w_7(\infty)=-9/20$
at which $t=-1073152081/3024000000$, a CM point 7-isogenous with
$t=\infty$ and thus of discriminant $-3\cdot 7^2$.

For each of $l=5,7,13$ the two fixed points of~$w_l$ on~$\XXo(l)$
are rational and yields two new CM points of discriminants $-cl$
for some factors~$c$ of~$24$.  For $\XXo(5)$ these
fixed points are $x_5=-3/5$ and $x_5=7/30$, at which
$t=2312/125$ and $t=5776/3375$ respectively; these CM points
have discriminants $-40$, $-120$ by the supersingular test.
For $\XXo(7)$ we find $x_7=2$ and $x_7=-29/10$, and thus
$t=-169/27$, $t=-701784/15625$ of discriminants $-84$, $-168$
divisible by~$7$.  For $\XXo(13)$ the fixed points $x_{13}=9$,
$x_{13}=-4$ yield $t=6877/15625$ and $t=27008742384/27680640625$,
with discriminants $-52=4\cdot13$ and $-312=24\cdot13$.

Each of these new CM points admits an $l$-isogeny to itself.
By solving the equation $t(x_l)=t(w_l(x_l))$ we find the remaining
such points; those not accounted for by fixed points of $w_l$ admit
two self-isogenies of degree~$l$, and correspond to a quadratic
pair of $x_l$ values over $\Q(t)$.  As it happens all the $t$\/'s
thus obtained are rational with the exception of a quadratic pair
coming from the quartic
$167x_{13}^4-60x_{13}^3+12138x_{13}^2-1980x_{13}+221607=0$.
Those points are: from $\XXo(5)$, the known $t=1$, $t=-169/25$,
and the new $t=-1377/1024$, $t=3211/1024$ of discriminants
$-51$, $-19$; from $\XXo(7)$, the CM points $t=0$, $152881/138240$,
$3211/1024$, $2312/125$, $6877/15625$ seen already, but also
$t=13689/15625$ of discriminant $-132$; and from $\XX_0(13)$,
seven of the CM points already known and also the two new
values $t=21250987/16000000$, $15545888/20796875$
of discriminants $-43$, $-88$.

\subsection{Numerical computation of CM points on~$\XX(1)$}

If we could obtain equations for the modular cover of~$\XX(1)$
by the elliptic curve $\XX(11)$, $\XX(17)$ or $\XX(19)$ we could
similarly find a few more rational CM~points on~$\XX(1)$.  But
we do not know how to find these covers, let alone the cover
$\XX(l)$ for $l$\/ large enough to get at the rational CM
point of discriminant $-163$; moreover, some applications
may require irrational CM~points of even higher discriminants.
We thus want a uniform way of computing the CM points of any
given discriminant as an algebraic irrationality.  We come close
to this by finding these points and their algebraic conjugates
as real (or, in the irrational case, complex) numbers to high
precision, and then using continued fractions to recognize
their elementary symmetric functions as rational numbers.
We say that this ``comes close'' to solving the problem because,
unlike the case of the classical modular functions such as~$j$,
we do not know a priori how much precision is required, since
the CM~values are generally not integers, nor is an effective
bound known on their height.  However, even when we cannot
prove that our results are correct using an isogeny of low
degree, we are quite confident that the rational numbers
we exist are correct because they not only match their
numerical approximations to many digits but also pass all
the supersingularity tests we tried as well as the condition
that differences between pairs of CM~values 
are products of small primes as in~\cite{GZ}.

To do this we must be able to compute numerically the
rational function~$t:\H/\Gamma^*(1){\stackrel\sim\ra}\Pr^1$.
Equivalently, we need to associate to each $t\in\Pr^1$ a representative
of its corresponding $\Gamma^*(1)$-orbit in~$\H$.  We noted already
that this is done, up to a fractional linear transformation over~$\C$,
by the quotient of two hypergeometric functions in~$t$.  To fix the
transformation we need images of three points, and we naturally
choose the elliptic points $t=0,1,\infty$.  These go to
fixed points of $s_2,s_4,s_6\in\Gamma^*(1)$, and to find those
fixed points we need an explicit action of~$\Gamma^*(1)$ on~$\H$.

To obtain such an action we must imbed that group
into $\Aut(\H) = \PSL_2(\R)$.  Equivalently, we
must choose an identification of $\A\otimes\R$ with the algebra
$M_2(\R)$ of $2\times 2$ real matrices.
Having done this, to obtain the action of some
$g\in\Gamma^*(1)\subset\A^*/\Q^*$ on~$\H$ we will choose a
representative of~$g$ in~$\A^*$, identify this representative
with an invertible matrix $({a\;b\atop c\;d})$ of positive determinant,
and let $g$ act on $z\in\H$ by $z \mapsto (az+b) / (cz+d)$.
Identifying $\A\otimes\R$ with $M_2(\R)$ is in turn tantamount
to solving (\ref{bcrel}) in $M_2(\R)$.

We choose the following solution:
\be
b:= \left(\begin{array}{cc}\sqrt2 & 0 \\ 0 & -\sqrt2 \end{array}\right),
\qquad
c:= \left(\begin{array}{cc}0 & \sqrt3 \\ -\sqrt3 & 0 \end{array}\right).
\label{b,c}
\ee
The elliptic points are then the $\Gamma^*(1)$ orbits of the fixed
points in the upper half-plane of $s_2,s_4,s_6$, that is, of
\be
P_2 := (1+\sqrt2)i,\quad P_4 := \frac{1+\sqrt2}{\sqrt3} (-1+\sqrt2\,i),
\quad P_6 := i.
\label{P246}
\ee
Thus for $|t|<1$ the point on~$\H/\Gamma^*(1)$ which maps to~$t$
is the $\Gamma^*(1)$ orbit of $z$ near $P_2$ such that
\be
(z-P_2) / (z-\bar P_2) = F_1(t) / F_2(t)
\label{param}
\ee
for some solutions $F_1,F_2$ of the hypergeometric
equation~(\ref{hypereq}).  Since the fractional linear
transformation $z\mapsto (z-P_2)/(z-\bar P_2)$
takes the hyperbolic lines $\overline{P_2 P_4}$ and
$\overline{P_2 P_6}$ to straight lines through the origin,
$F_2$ must be a power series in~$t$, and $F_1$ is such a
power series multiplied by $\sqrt{t}$; that is,
\be
(z-P_2) / (z-\bar P_2)
=
C t^{1/2}
 F\left(\frac{13}{24},\frac{17}{24},\frac32,t\right)
 \left/
 F\left(\frac1{24},\frac5{24},\frac12,t)\right)
 \right. .
\label{param1}
\ee
for some nonzero constant~$C$.  We evaluate $C$\/ by taking
$t=1$ in~(\ref{param1}).  Then $z=P_4$, which determines
the left-hand side, while the identity~\cite[9.122]{GR}
\be
F(a,b;c;1) = \frac{\Gamma(c) \Gamma(c-a-b)} {\Gamma(c-a) \Gamma(c-b)}
\label{hypid}
\ee
gives us the coefficient of $C$\/ in the right-hand side in terms
of gamma functions.  We find $C=(.314837\ldots)i/(2.472571\ldots)
= (.128545\ldots) i$.  Likewise we obtain convergent power series for
computing $z$ in neighborhoods of $t=1$ and $t=\infty$.

Now let $D$\/ be the discriminant of an order $O_D$ in a quadratic
imaginary field $\Q(\sqrt{D}\!)$ such that $O_D$ has a maximal embedding
in~$\O$ (i.e.\ an embedding such that $O_D = (O_D\otimes\Q)\cap\O$)
and the embedding is unique up to conjugation in $\Gamma^*(1)$.
Then there is a unique, and therefore rational, CM point on~$\XX(1)$
of discriminant~$D$.  Being rational, the point is real, and thus
can be found on one of the three hyperbolic line segments
$\overline{P_2 P_4}$, $\overline{P_2 P_6}$, $\overline{P_4 P_6}$.
It is thus the fixed point of a positive integer combination, with
coprime coefficients, of two of the elliptic elements 
$s_2 = bc+2c$, $s_4 = (2+b)(1+c)/2$, $s_6=(3+c)/2$ with
fixed points $P_2,P_4,P_6$.  In each case a short search
finds the appropriate linear combination and thus the
fixed point~$z$.  Using (\ref{param1}) or the analogous
formulas near $t=1$, $t=\infty$ we then solve for $t$
as a real number with sufficient accuracy (60 decimals
was more than enough) to recover it as a rational number
from its continued-fraction expansion.

\subsection{Tables of rational CM points on~$\XX(1)$}

There are 27 rational CM points on $\XX(1)$.  We write
the discriminant~$D$\/ of each of them as $-D_0 D_1$
where $D_0|24$ and $D_1$ is coprime to~$6$.
In Table~1 we give, for each $|D|=D_0 D_1$, the integers $A,B$
with $B\geqs0$ such that $(A:B)$ is the $t$-coordinate of
a CM~point of discriminant~$D$.  In the last column of
this table we indicate whether the point was obtained
algebraically (via an isogeny of degree 5, 7, or~13)
and thus proved correct, or only computed numerically.  The
CM~points are listed in order of increasing height $\max(|A|,B)$.

\begin{figure}[t]
\centerline{\large Table 1}

\vspace*{2ex}

\def\m{\phantom-}

\setlength{\tabcolsep}{4pt}
\centerline{
\begin{tabular}{c|c|c|c|c|c}
$|D|$&$D_0$&$D_1$&    $A$        &    $B$    & proved? \\ \hline
  3  &  3  &  1  &     1         &     0     &  Y  \\
  4  &  4  &  1  &     1         &     1     &  Y  \\
 24  & 24  &  1  &     0         &     1     &  Y  \\
 84  & 12  &  7  &  $-169\m$     &    27     &  Y  \\
 40  &  8  &  5  &    2312       &    25     &  Y  \\
 51  &  3  & 17  &  $-1377\m$    &   1024    &  Y  \\
 19  &  1  & 19  &    3211       &   1024    &  Y  \\
120  & 24  &  5  &    5776       &   3375    &  Y  \\
 52  &  4  & 13  &    6877       &   15625   &  Y  \\
132  & 12  & 11  &   13689       &   15625   &  Y  \\
 75  &  3  &$5^2$&   152881      &  138240   &  Y  \\
168  & 24  &  7  & $-701784\m$   &   15625   &  Y  \\
 43  &  1  &  43 &  21250987     & 16000000  &  Y  \\
228  & 12  &  19 &  66863329     & 11390625  &  N  \\
 88  &  8  &  11 &  15545888     & 20796875  &  Y  \\
123  &  3  &  41 &$-296900721\m$ & 16000000  &  N  \\
100  &  4  &$5^2$&  421850521    &  1771561  &  Y  \\
147  &  3  &$7^2$&$-1073152081\m$&3024000000 &  Y  \\
312  & 24  &  13 & 27008742384   &27680640625&  Y  \\
 67  &  1  &  67 & 77903700667   &1024000000 &  N  \\
148  &  4  &  37 & 69630712957       & 377149515625 &  N  \\
372  & 12  &  31 & $-455413074649\m$ & 747377296875 &  N  \\
408  & 24  &  17 &$-32408609436736\m$&  55962140625 &  N  \\
267  &  3  &  89 &$-5766681714488721\m$&   1814078464000000  & N \\
232  &  8  &  29 & 66432278483452232   &  56413239012828125  & N \\
708  & 12  &  59 & 71475755554842930369&224337327397603890625& N \\
163  &  1  & 163 &699690239451360705067&684178814003344000000& N
\end{tabular}
}
\setlength{\tabcolsep}{1.4pt}
\end{figure}

In Table~2 we give, for each except the first three cases,
the factorizations of $|A|,B,|C|$ where $C=A-B$, and also
the associated ``$\!ABC$\/~ratio''~\cite{ABC->Mord} defined by
$r=\log N(ABC)/\log\max(|A|,B,|C|)$.  As expected, the
$A,B,C$\/~values are ``almost'' perfect squares, sixth powers,
and fourth powers respectively: a prime at which at which the
valuation of $A,B,C$\/ is not divisible by 2,\,6,\,4 resp.\
is either 2, 3, or the unique prime in $D_1$.  When $D_1>1$
its unique prime factor is listed at the end of the $|A|$,
$B$, or $|C|$ factorization in which it appears; otherwise
the prime factors are listed in increasing order.

\begin{figure}[t]
\centerline{\large Table 2}

\vspace*{2ex}
\setlength{\tabcolsep}{4pt}
\centerline{
$
\begin{array}{c|c|c|c|c|c|c}
|D|& D_0 & D_1 &    |A|      &     B     &    |C|   & r \\ \hline
 84& 12  &  7  &   13^2      &    3^3    &  2^2 7^2 & 1.19410 \\
 40&  8  &  5  &  2^3 17^2   &    5^3    &    3^7   & 0.80487 \\
 51&  3  & 17  &   3^4 17    &  2^{10}   &    7^4   & 0.84419 \\
 19&  1  & 19  &  13^2 19    &  2^{10}   &    3^7   & 0.90424 \\
120& 24  &  5  &  2^4 19^2   &  3^3 5^3  &    7^4   & 0.95729 \\
 52&  4  & 13  &  23^2 13    &    5^6    &  2^2 3^7 & 1.00276 \\
132& 12  & 11  &  3^4 13^2   &    5^6    &  2^4 11^2& 0.87817 \\
 75&  3  & 5^2 &  17^2 23^2  &2^{10}3^3 5&    11^4  & 0.98579 \\
168& 24  &  7  & 2^3 3^5 19^2&    5^6    & 11^4 7^2 & 0.79278 \\
 88&  8  &  11 &2^5 17^2 41^2&  5^6 11^3 &  3^7 7^4 & 0.86307 \\
 43&  1  &  43 & 19^2 37^2 43& 2^{10}5^6 &  3^7 7^4 & 0.92839 \\
228& 12  &  19 &13^2 17^2 37^2  &   3^6 5^6    &2^6 7^4 19^2 & 0.96018\\
123&  3  &  41 &3^4 13^2 23^2 41&  2^{10}5^6   &   7^4 19^4  & 0.90513\\
100&  4  & 5^2 &19^2 23^2 47^2  &    11^6      &2^4 3^7 7^4 5& 0.88998\\
147&  3  & 7^2 &17^2 41^2 47^2 &2^{10}3^3 5^6 7&  11^4 23^4  & 0.96132\\
312& 24  &  13 &2^4 3^5 17^2 43^2 13& 5^6 11^6 &  7^4 23^4  & 0.83432\\
 67&  1  &  67 &  13^2 43^2 61^2 67 &2^{16}5^6 &3^7 7^4 11^4 & 0.89267\\
148&  4  &  37 &  13^2 47^2 71^2 37 & 5^6 17^6 &
                2^2 3^7 7^4 11^4 & 0.94008 \\
372& 12  &  31 &13^2 23^2 37^2 61^2 &3^3 5^6 11^6 &
                2^2 7^4 19^4 31^2 & 0.99029 \\
408& 24  &  17 &2^6 13^2 19^2 43^2 67^2 &3^6 5^6 17^3 &
                7^4 11^4 31^4 & 0.88352 \\
267&  3  &  89 &3^6 13^2 17^2 19^2 71^2 89 & 2^{16}5^6 11^6 &
                7^4 31^4 43^4  &   0.87610 \\
232&  8  &  29 &2^3 13^2 17^2 41^2 89^2 113^2 &5^6 23^6 29^3&
                3^7 7^4 11^4 19^4  &   0.91700 \\
708& 12  &  59 &3^4 13^2 19^2 23^2 37^2 41^2 109^2 &  5^6 17^6 29^6 &
                2^8 7^4 11^4 47^4 59^2  &   0.91518 \\
163&  1  & 163 &13^2 67^2 109^2 139^2 157^2 163 & 2^{10} 5^6 11^6 17^6 &
                3^{11} 7^4 19^4 23^4 & 0.90013
\end{array}
$
}
\setlength{\tabcolsep}{1.4pt}
\end{figure}

In the factorization of the difference between the last two $t=A/B$\/
values in this table, the primes not accounted for by common factors
in the last two rows of the table are 79, 127, 271, 907, 2287, 2971,
3547, each occurring once.

\section{The case $\Sigma = \{2,5\}$}

\subsection{The quaternion algebra and the curves $\Xone$, $\XX(1)$}

For this section we let $\A$ be the quaternion algebra ramified
at $\{2,5\}$.  This time $\A$ is generated over~$\Q$ by elements
$b,e$ satisfying
\be
b^2+2 = e^2-5 = be+eb = 0,
\label{berel}
\ee
and the conjugate and norm of an element
$\alpha = \alpha_1 + \alpha_2 b + \alpha_3 e + \alpha_4 be \in \A$
are
\be
\bar\alpha = \alpha_1 - \alpha_2 b - \alpha_3 e - \alpha_4 be,
\quad
\alpha\bar\alpha = \bar\alpha\alpha =
\alpha_1^2 + 2\alpha_2^2 - 5\alpha_3^2 - 10\alpha_4^2.
\ee
The elements $b$ and $(1+e)/2$ generate a maximal order, which
we use for~$\O$.

By (\ref{Area*}), the curve $\XX(1)$ has hyperbolic area $1/6$.
Since the algebra $\A$ is not among the nineteen algebras listed
in~\cite{Takeuchi} that produce arithmetic triangle groups, $\XX(1)$
must have at least four elliptic points.  On the other hand,
by~(\ref{genus}) a curve of area as small as $1/6$ cannot have
more than four elliptic points, and if it has exactly four then
their orders must be $2,2,2,3$.
Indeed we find in $\Gamma^*(1)$ the elements of finite order
\be
s_2 = [b], s'_2 = [2e + 5b - be], s''_2 = [5b - be], s_3 = [2b - e - 1]
\label{s2223}
\ee
[NB $2e+5b-be, 5b-be, 2b-e-1\in 2\O$] of orders $2,2,2,3$ with
$s_2^{\phantom0} s'_2 s''_2 s_3^{\phantom0} = 1$.
As in the case of the $G_{2,4,6}$ we
conclude that here $\Gamma^*(1)$ has the presentation
\be
\langle s_2^{\phantom0},s'_2,s''_2,s_3^{\phantom0} |
s_2^2 = {s'_2}^2 = {s''_2}^2 = s_3^3
= s_2^{\phantom0} s'_2 s''_2 s_3^{\phantom0} = 1
\rangle
.
\label{G2223}
\ee
Of the four generators only $s_3$ is in $\Gamma(1)$; thus the
$(\Z/2)^2$ cover $\Xone/\XX(1)$ is ramified at the elliptic points
of order~2.   Therefore $\Xone$ is a rational curve with four elliptic
points of order~3, and $\Gamma^(1)$ is generated by four 3-cycles
whose product is the identity, for example by $s_3$ and its
conjugates by $s_2,s'_2,s''_2$.  (The genus and number of elliptic
points of $\Xone,\XX(1)$, but not the generators of
$\Gamma(1),\Gamma^*(1)$, are already tabulated in \cite[Ch.IV:2]{V}.)

\subsection{Shimura modular curves $\XXo(l)$, in particular $\XXo(3)$}

The elliptic elements $s_3^{\phantom0},s^{\phantom0}_2,s'_2,s''_2$ have
discriminants $-3,-8,-20,-40$.  Thus the curve $\XXo(l)$ has genus
\be
g(\XXo(l)) = \frac1{12} \left[ l
- 4 \left(\frac{-3}{l}\right)
- 3 \left(\frac{-2}{l}\right)
- 3 \left(\frac{-5}{l}\right)
- 3 \left(\frac{-10}{l}\right)
\right] .
\label{genus(10)}
\ee
Again we tabulate this for $l<50$:

\vspace*{2ex}

\setlength{\tabcolsep}{4pt}
\centerline{
\begin{tabular}{c|ccccccccccccc}
$l$ &
3 & 7 & 11 & 13 & 17 & 19 & 23 & 29 & 31 & 37 & 41 & 43 & 47 \\ \hline
$g(\XXo(l))$ & 
0 & 0 &  1 &  1 &  2 &  1 &  2 &  3 &  3 &  3 &  3 &  3 &  4
\end{tabular}
}
\setlength{\tabcolsep}{1.4pt}

\vspace*{2ex}

Since $g(\XXo(l))\geqs(l-13)/12$, the cases $l=3,7$ of genus~0
occurring in this table are the only ones.
We next find an explicit rational functions of degree~4 on~$\Pr^1$
that realizes the cover $\XXo(3)/\XXo(1)$,
and determine the involution~$w_3$.

The curve $\XXo(3)$ is a degree-4 cover of $\XX(1)$ with Galois
group $\PGL_2(\F_3)$ and cycle structures 31, 211, 211, 22 over
the elliptic points $P_3^{\phantom0},P^{\phantom0}_2,P'_2,P''_2$.
Thus there are coordinates $\tau,x$ on $\XX(1)$, $\XXo(3)$ such
that $\tau(x)=(x^2-c)^2/(x-1)^3$ for some~$c$.
To determine the parameter~$c$, we use the fact that $w_3$
fixes the simple pole $x=\infty$ and takes each simple preimage of the
211~points $P_2^{\phantom0},P'_2$ to the other simple preimage
of the same point.  That is,
\be
(x^2-c)^{-1} (x-1)^4 \frac{dx}{dt} = x^2 - 4x + 3c
\ee
must have distinct roots $x_i$ ($i=1,2$) that yield quadratic
polynomials
\be
\frac{(x-1)^3 (\tau(x)-\tau(x_i))} {(x-x_i)^2}
\ee
with the same $x$ coefficient.  We find that this happens
only for $c=-5/3$, i.e.\ that $\tau=(3x^2+5)^2/9(x-1)^3$.
For future use it will prove convenient to use
\be
t = \frac{6^3}{9\tau+8} = \frac{(6x-6)^3}{(x+1)^2(9x^2-10x+17)},
\label{x10:3}
\ee
with $w_3(x)=\frac{10}{9}-x$.  [Smaller coefficients can be
obtained by letting $x=1+2/x'$, $\tau=2t'/9$, when
$t'=(2x'^2+3x'+3)^2/x'$ and $w_3(x')=-9x'/(4x'+9)$.  But
our choice of $x$ will simplify the computation of the
Schwarzian equation, while the choice of~$t$\/ will turn out
to be the correct one 3-adically.]  The elliptic points are then
$P_3:t=0$, $P''_2:t=27$, and $P_2,P'_2:t=\infty,2$.  In fact
the information so far does not exclude the possibility that
the pole of~$t$ might be at $P'_2$ instead of $P_2$; that
in fact $t(P_2)=\infty,t(P'_2)=2$ and not the other way around
can be seen from the order of the elliptic points on the real
locus of~$\XX(1)$, or (once we compute the Schwarzian equation)
checked using the supersingular test.

\subsection{CM points on~$\XX(1)$ via $\XXo(3)$ and $w_3$}

{}From $w_3$ we obtain five further CM points.  Three
of these are 3-isogenous to known elliptic points: $w_3$
takes the triple zero $x=1$ of~$t$\/ to $x=1/9$, which gives us
$t=-192/25$, the point 3-isogenous to $P_3$ with discriminant~$-27$;
likewise $w_3$ takes the double root $x=5$ and double pole $x=-1$
of~$t-2$ to $x=-35/9,19/9$ and thus to $t=-2662/169$ and $t=125/147$,
the points 3-isogenous to $t=2$ and $t=\infty$ and thus (once these
points are identified with $P'_2$ and $P_2$) of discriminants
$-180$ and $-72$.  One new CM point comes from the other fixed point
$x=5/9$ of $w_3$, which yields $t=-27/49$ of discriminant~$-120$.
Finally the remaining solutions of $t(x)=t(w_3(x))$ are the roots
of $9x^2-10x+65$; the resulting CM point $t=64/7$, with two
3-isogenies to itself, turns out to have discriminant~$-35$.

\subsection{The Schwarzian equation on~$\XX(1)$}

We can take the Schwarzian equation on $\XX(1)$ to be of the form
\be
t (t-2) (t-27) f'' + (A t^2 + B t + C) f' + (D t + E) = 0.
\label{schraw}
\ee
The coefficients $A,B,C,D$\/ are then forced by the indices
of the elliptic points.  Near $t=0$, the solutions of~(\ref{schraw})
must be generated by functions with leading terms $1$ and~$t^{1/3}$;
near $t=2$ ($t=27$), by functions with leading terms 1 and
$(t-2)^{1/2}$ (resp.\ $(t-27)^{1/2}$); and at infinity, by
functions with leading terms $t^{-e}$ and $t^{-e-1/2}$ for
some~$e$.  The conditions at the three finite singular points
$t=0,2,27$ determine the value of the $f'$ coefficient at those
points, and thus yield $A,B,C$, which turn out to be $5/3,-203/6,36$.
Then $e,e+1/2$ must be roots of an ``indicial equation''
$e^2-2e/3+D=0$, so $e=1/12$ and $D=7/144$.

Thus (\ref{schraw}) becomes
\be
t (t-2) (t-27) f'' + \frac{10t^2-203t+216}{6} f' + (\frac{7t}{144} + E)
= 0.
\label{schw_E}
\ee
To determine the ``accessory parameter'' $E$, we again use
the cover $\XXo(3)/\XX(1)$ and the involution $w_3$.  A Schwarzian
equation for $\XXo(3)$ is obtained by substituting
$t=(6x-6)^3/(x+1)^2(9x^2-10x+17)$ in~(\ref{schw_E}).  The resulting
equation will not yet display the $w_3$ symmetry, because it will
have a spurious singular point at the double pole $x=-1$ of~$t(x)$.
To remove this singularity we consider not $f(t(x))$ but
\be
g(x) := (x+1)^{-1/6} f(t(x)).
\label{f2g}
\ee
The factor $(x+1)^{-1/6}$ is also singular at $x=\infty$, but that
is already an elliptic point of $\XXo(3)$ and a fixed point of $w_3$.
Let $x=u+5/9$, so $w_3$ is simply $u\leftrightarrow -u$.  Then
we find that the differential equation satisfied by $g$ is
$$
4 (81u^2+20) (81u^2+128)^2 g''
+ 108 u (81u^2+128) (405u^2+424) g'
\qquad
$$ \vspace*{-6ex} \be \label{schw_g} \ee \vspace*{-4ex} $$
\quad + (3^{11}u^4-163296u^2+170496+72(18E+7)(9u-4)(81u^2+128)) g
= 0.
$$
Clearly this has the desired symmetry if and only if $18E+7=0$.
Thus the Schwarzian equation is
\be
t (t-2) (t-27) f'' + \frac{10t^2-203t+216}{6} f' +
(\frac{7t}{144} - \frac{7}{18}) = 0.
\label{schw:10}
\ee

\subsection{Numerical computation of CM points on~$\XX(1)$}

We can now expand a basis of solutions of (\ref{schw:10}) in power
series about each singular point $t=0,2,27,\infty$ (using inverse
powers of $t-\frac{27}{2}$ for the expansion about $\infty$ to
assure convergence for real $t\notin[0,27]$).  As with the $\Sigma=\{2,3\}$ case
we need to identify $\A\otimes\R$ with $M_2(\R)$, and use the
solution
\be
b:= \left(\begin{array}{cc}0 & \sqrt2 \\ -\sqrt2 & 0 \end{array}\right),
\qquad
e:= \left(\begin{array}{cc}\sqrt5 & 0 \\ 0 & -\sqrt5 \end{array}\right).
\label{b,e}
\ee
of~(\ref{berel}), analogous to (\ref{bcrel}).  We want to proceed
as we did for $\Sigma=\{2,3\}$, but there is still one obstacle to
computing, for given $t_0\in\R$, the point on the hyperbolic
quadrilateral formed by the fixed points of
$s_2^{\phantom0},s'_2,s''_2,s_3^{\phantom0}$
at which $t=t_0$.  In the $\Sigma=\{2,3\}$ case, the solutions
of the Schwarzian equation were combinations of hypergeometric
functions, whose value at~$1$ is known.  This let us determine two
solutions whose ratio gives the desired map to~$\H$.  But here
$\Gamma^*(1)$ is not a triangle group, so our basic solutions
of~(\ref{schw:10} are more complicated power series and we
do not know a priori their values at the neighboring singular points.
In general this obstacle can be overcome by noting that for each
nonsingular $t_0\in\R$ its image in~$\H$ can be computed from the
power-series expansions about either of its neighbors and using
the condition that the two computations agree for several
choices of $t_0$ to determine the maps to~$\H$.
In our case we instead removed the obstacle using the non-elliptic
CM points computed in the previous section.  For example,
we used the fact that $t_0=125/147$ is the CM point
of discriminant~$72$, and thus maps to the unique fixed point in~$\H$
of $(9b+4e-be)/2$, to determine the correct ratio of power series
about $t=0$ and $t=2$.  Two or three such points suffice to determine
the four ratios needed to compute our map $\R\ra\H$ to arbitrary
accuracy; since we actually had five non-elliptic CM points,
we used the extra points for consistency checks, and then used
the resulting formulas to numerically compute the $t$-coordinates
of the remaining CM points.

There are 21 rational CM points on $\XX(1)$.  We write
the discriminant~$D$\/ of each of them as $-D_0 D_1$
where $D_0|40$ and $D_1$ is coprime to~10.  Table~3 is
organized in the same way as Table~1: we give, for each $|D|=D_0 D_1$,
the integers $A,B$ with $B\geqs0$ such that $(A:B)$ is the
$t$-coordinate of a CM~point of discriminant~$D$.
The last column identifies with a ``Y'' the nine points
obtained algebraically from the computation of $\XXo(3)$ and $w_3$.
Some but not all of the remaining twelve points would move from
``N'' to ``Y'' if we also had the equations for the degree-8
map $\XXo(7)\ra\XX(1)$ and the involution $w_7$ on~$\XXo(7)$.

\begin{figure}[t]
\setlength{\tabcolsep}{4pt}
\centerline{\large Table 3}

\vspace*{2ex}

\def\m{\phantom-}

\centerline{
\begin{tabular}{c|c|c|c|c|c}
$|D|$&$D_0$&$D_1$&      $A$      &    $B$    & proved? \\ \hline
  3  &  1  &  1  &       0       &     1     &  Y  \\
  8  &  8  &  1  &       1       &     0     &  Y  \\
 20  & 20  &  1  &       2       &     1     &  Y  \\
 40  & 40  &  1  &      27       &     1     &  Y  \\
 52  &  4  & 13  &    $-54\m$    &    25     &  N  \\
120  & 40  &  3  &    $-27\m$    &    49     &  Y  \\
 35  &  5  &  7  &      64       &     7     &  Y  \\
 27  &  1  &$3^3$&   $-192\m$    &    25     &  Y  \\
 72  &  8  &$3^2$&     125       &    147    &  Y  \\
 43  &  1  & 43  &    1728       &   1225    &  N  \\
180  & 20  &$3^2$&  $-2662\m$    &    169    &  Y  \\
 88  &  8  & 11  &    3375       &    98     &  N  \\
115  &  5  & 23  &   13824       &   3887    &  N  \\
280  & 40  &  7  &   35937       &   7406    &  N  \\
 67  &  1  & 67  & $-216000\m$   &   8281    &  N  \\
148  &  4  & 37  &    71874      &  207025   &  N  \\
340  & 20  & 17  &    657018     &  41209    &  N  \\
520  & 40  & 13  &    658503     & 11257064  &  N  \\
232  &  8  & 29  &   176558481   &  2592100  &  N  \\
760  & 40  & 19  &  13772224773  & 237375649 &  N  \\
163  &  1  &163  &$-2299968000\m$&6692712481 &  N  \\
\end{tabular}
}
\setlength{\tabcolsep}{1.4pt}
\end{figure}

It will be seen that the factor $3^3$ in our normalization
(\ref{x10:3}) of~$t$\/ was needed\footnote{
  On the other hand the factor $2^3$ in~(\ref{x10:3}) was
  a matter of convenience, to make the four elliptic points integral.
  }
to make $t$ a good coordinate 3-adically:
$3$ splits in the CM~field iff $t$ is not a multiple of~$3$.

In Table~4 we give the factorizations of $|A|,B,|A-2B|,|A-27B|$;
as expected, $|A|$ is always ``almost'' a perfect cube, and
$B,|A-2B|,|A-27B|$ ``almost'' a perfect square, any exceptional
primes other than 2 or~5 being the unique prime in~$D_1$, which
if it occurs is listed at the end of its respective factorization.

\begin{figure}[t]
\centerline{\large Table 4}

\vspace*{2ex}

\setlength{\tabcolsep}{4pt}
\centerline{
\def\M{\!\cdot\!}
$
\begin{array}{c|c|c|c|c|c|c}
|D|& D_0 & D_1 &  |A|    &   B   & |A-2B| & |A-27B| \\ \hline
  3&  1  &  1  &   0     &   1   &   2    &   3^3   \\
  8&  8  &  1  &   1     &   0   &   1    &    1    \\
 20& 20  &  1  &   2     &   1   &   0    &   5^2   \\
 40& 40  &  1  &  3^3    &   1   &  5^2   &    0    \\
 52&  4  & 13  & 2\M3^3  &  5^2  & 2^3 13 &   3^6   \\
120& 40  &  3  &  3^3    &  7^2  &  5^3   &2\M3^3 5^2\\
 35&  5  &  7  &  2^6    &   7   &2\M5^2  &   5^3     \\
 27&  1  & 3^3 & 2^6 3   &  5^2  &2\M11^2 & 17^2 3    \\
 72&  8  & 3^2 &  5^3    & 7^2 3 &  13^2  & 2^2 31^2  \\
 43&  1  & 43  & 2^6 3^3 &5^2 7^2&2\M19^2 &  3^6 43   \\
180& 20  & 3^2 & 2\M11^3 & 13^2  &2^3 5^3 3& 5^2 17^2  \\
 88&  8  & 11  & 3^3 5^3 &2\M7^2 &17^2 11 &  3^6      \\
115&  5  & 23  & 2^9 3^3 &13^2 23&2\M5^2 11^2&3^6 5^3\\
280& 40  &  7  &3^3 11^3 &2\M23^2 7&5^3 13^2 &3^8 5^2\\
 67&  1  & 67  &2^6 3^3 5^3& 7^2 13^2 &2\M11^2 31^2 &3^8 67\\
148&  4  & 37  &2\M3^3 11^3&5^2 7^2 13^2&2^5 17^2 37&3^8 29^2\\
340& 20  & 17  &2\M3^3 23^3&  7^2 29^2  &2^3 5^2 13^2 7&3^6 5^4\\
520& 40  & 13  &3^3 29^3 & 2^3 7^2 47^2 13 &5^4 11^2 17^2 &3^8 5^2 43\\
232&  8  & 29  &3^3 11^3 17^3&2^2 5^2 7^2 23^2&13^2 19^2 53^2&3^6 71^2 29\\
760& 40  & 19  &  3^3 17^3 47^3  &   7^2 31^2 71^2    &
   5^2 11^2 13^2 37^2 19& 2\M3^8 5^3 67^2 \\
163&  1  &163  &2^9 3^3 5^3 11^3 & 7^2 13^2 29^2 31^2 & 
   2\M19^2 59^2 79^2 & 3^6 17^2 73^2 163 \\
\end{array}
$
}
\setlength{\tabcolsep}{1.4pt}
\end{figure}

\section{Further examples and problems}

Our treatment here is briefer because most of the ideas and
methods of the previous sections apply here with little change.
Thus we only describe new features that did not arise for the
algebras ramified at $\{2,3\}$ and $\{2,5\}$, and exhibit
the final results of our computations of modular curves and CM~points.

\subsection{The case $\Sigma = \{2,7\}$}

We generate $\A$ by elements $b,g$ with
\be
b^2+2 = g^2-7 = bg+gb = 0,
\label{bgrel}
\ee
and a maximal order~$\O$ by $\Z[b,g]$ together with
$(1+b+g)/2$ (and $b(1+g)/2$).  
By (\ref{Area*}), the curve $\XX(1)$ has hyperbolic area $1/4$.
Since $\Gamma^*(1)$ is not a triangle group (again by~\cite{Takeuchi}),
we again conclude by~(\ref{genus}) that $\XX(1)$ has exactly four
elliptic points, this time of orders $2,2,2,4$.  We find in
$\Gamma^*(1)$ the elements of finite order
\be
s_2 = [b], s'_2 = [7b - 2g - bg], s''_2 = [7b + 2g - bg],
s_4 = [1 + 2b + g]
\label{s2224}
\ee
[NB $7b \pm 2g - bg\in 2\O$] of orders $2,2,2,4$ with
$s_2^{\phantom0} s'_2 s''_2 s_4^{\phantom0} = 1$,
and conclude that $s_2^{\phantom0},s'_2,s''_2,s_4^{\phantom0}$
generate $\Gamma^*(1)$ with relations determined by
$s_2^2 = {s'_2}^2 = {s''_2}^2 = s_4^4 = s_2 s'_2 s''_2 s_4 = 1$.
None of these is in $\Gamma^*(1)$: the representatives
$b,1+2b+g$ of $s_2,s_4$ have norm~$2$, while $s'_2,s''_2$
have representatives $(7b \pm 2g - bg)/2$ of norm~$14$.
The discriminants of $s_4,s_2,s'_2,s''_2$ are $-4,-8,-56,-56$;
note that $-56$ is {\em not}\/ among the ``idoneal'' discriminants
(discriminants of imaginary quadratic fields with class group
$(\Z/2)^r$), and thus that the elliptic fixed points $P'_2,P''_2$
of $s'_2,s''_2$ are quadratic conjugates on~$\XX(1)$.

Again we use the involution $w_3$ on the modular curve $\XXo(3)$
to simultaneously determine the relative position of the elliptic
points $P_4,P_2,P'_2,P''_2$ on $\XX(1)$ and the modular cover
$\XXo(3)\ra\XX(1)$, and then to obtain a Schwarzian equation
on~$\XX(1)$.  Clearly $P_4$ is completely ramified in $\XXo(3)$.
Since $-8$ and $-56$ are quadratic residues of~$3$, each of
$P_2,P'_2,P''_2$ has ramification type 211.  Thus $\XXo(3)$
is a rational curve with six elliptic points all of index~2,
and we may choose coordinates $t,x$ on $\XX(1),\XXo(3)$ such
that $t(P_4)=\infty$, $t(P_2)=0$, and $x=\infty$, $x=0$ at
the quadruple pole and double zero respectively of~$t$.

We next determine the action of $w_3$ on the elliptic points
of~$\XXo(3)$.  Necessarily the simple preimages of $P_2$ parametrize
two 3-isogenies from $P_2$ to itself.  On the other hand the simple
preimages of $P'_2$ parametrize two 3-isogenies from that point to
$P''_2$ and vice versa, because the squares of the primes above~$3$
in $\Q(\sqrt{-14})$ are not principal.  Therefore $w_3$ exchanges
the simple preimages of $P_2$ but takes each of the two simple
points above~$P'_2$ to one above $P''_2$ and vice versa.

So again we have a one-parameter family of degree-4 functions
on~$\Pr^1$, and a single condition in the existence of the
involution $w_3$; but this time it turns out that there are
(up to scaling the coordinates $t,x$) two ways to satisfy
this condition:
\be
t = \frac13(x^4 + 4 x^3 + 6 x^2), \quad w_3(x) = \frac{1-x}{1+x}, \quad
P'_2,P''_2: t^2 - 3t + 3 = 0
\label{14:no}
\ee
and 
\be
t = \frac1{27}(x^4 + 2 x^3 + 9 x^2), \quad w_3(x) = \frac{5-2x}{2+x},
\quad P'_2,P''_2: 16 t^2 + 13 t + 8 = 0.
\label{14:yes}
\ee
How to choose the correct one?  We could consider the next modular
curve $\XXo(5)$ and its involution to obtain a new condition that
would be satisfied by only one of (\ref{14:no},\ref{14:yes}).
Fortunately we can circumvent this laborious calculation by
noting that the Fuchsian group associated with (\ref{14:no})
is commensurable with a triangle group, since its three elliptic
points of index~2 are the roots of $(1-t)^3=1$ and are thus
permuted by a 3-cycle that fixes the fourth elliptic point $t=\infty$.
The quotient by that 3-cycle is a curve parametrized by $(1-t)^3$
with elliptic points of order $2,3,12$ at $1,0,\infty$.  But
by~\cite{Takeuchi} there is no triangle group commensurable with an
arithmetic subgroup of $\A^*/\Q^*$; indeed we find there that
$G_{2,3,12}$ is associated with the quaternion algebra over
$\Q(\sqrt{3})$ ramified at the prime above~2 and at one of
the infinite places of that number field.\footnote{
  See \cite{Takeuchi}, table 3, row~IV.  In terms of that algebra
  $\A'$, the triangle group $G_{2,3,12}$ is $\Gamma^*(1)$; the index-3
  normal subgroup whose quotient curve is parametrized by the~$t$\/
  of~(\ref{14:no}) is the normalized in $\Gamma^*(1)$ of
  $\{ [a] \in \O^*/\{\pm1\} : a\equiv1\bmod I_2\}$; and the
  intersection of this group with $\Gamma^*_0(3)$ yields as quotient
  curve the $\Pr^1$ with coordinate~$x$.
  }
Therefore (\ref{14:yes}) is the correct choice.  Alternatively,
we could have noticed that since $\Xone$ is a $(\Z/2)^2$ cover
of $\XX(1)$ ramified at all four elliptic points, it has genus~1,
and then used the condition that this curve's Jacobian have
conductor~14 to exclude (\ref{14:no}).  The function field
of $\XX(1)$ is obtained by adjoining square roots of $c_0 t$
and $c_1(16 t^2 + 13 t + 8)$ for some $c_0,c_1$; for the
Jacobian to have the correct conductor we must have $c_0 c_1=1$
mod~squares.  The double cover of $\XXo(3)$ obtained by adjoining
$\sqrt{c_1(16 t^2 + 13 t + 8)}$ also has genus~1, and so must
have Jacobian of conductor at most~42; this happens only when
$c_1=-1$ mod~squares, the Jacobian being the elliptic curve
42-A3~(42C).  The curve $\Xone$ then has the equation
\be
y^2 = - 16 s^4 + 13 s^2 - 8
\qquad(t=-s^2),
\label{14:X1}
\ee
and its Jacobian is the elliptic curve 14-A2 (14D).  Kurihara had
already obtained in~\cite{Ku} an equation birational with~(\ref{14:X1}).
Let $\Gamma'_0(3^r)$ be the group intermediate between
$\Gamma_0(3^r)$ and $\Gamma^*_0(3^r)$ consisting of the elements
of norm~1 or~7 mod~${\Q^*}^2$.  Then the corresponding curves
$\calX'_0(3^r)$ ($r>0$) of genus $3^{r-1}+1$ are obtained from
$\XXo(3^r)$ by extracting a square root of $t(16t^2+13t+8)$,
and constitute an unramified tower of curves over the genus-2 curve
\be
\calX'_0(3): y^2 = 3(4x^6 + 12x^5 + 75x^4 + 50x^3 + 255x^2 - 288x + 648)
\label{14:X3'}
\ee
whose reductions are asymptotically optimal over $\F_{l^2}$
($l\neq2,3,7)$
% perhaps 2,7 don't need to be excluded?
with each step in the tower being a cyclic cubic extension.
(Of course when we consider only reductions to curves over $\F_{l^2}$
the factor of~3 in~(\ref{14:X3'}) may be suppressed.)

Using $w_3$ we may again find the coordinates of several non-elliptic
CM points: $t=4/3$ and $t=75/16$ of discriminants $-36$ and $-72$, i.e.\
the points {3-isogenous} to $P_4$ and $P_2$, other than $P_4,P_2$
themselves; $t=4/9$ and $t=200/9$ of discriminants $-84$ and $-168$,
coming from the fixed points $x=1$ and $x=-5$ of~$w_3$; and the
points $t=-1$, $t=-5$ of discriminants $-11$ and $-35$, coming from
the remaining solutions of $t(x)=t(w_3(x))$ and each with two
{3-isogenies} to itself.

Even once the relative position of the elliptic points is known,
the computation of the cover $\XXo(5)/\XX(1)$ is not a trivial
matter; I thank Peter M\"uller for performing this computation
using J.-C. Faugere's Gr\"obner basis package~GB.  It turns out
that there are eight PGL$_2(\F_5)$ covers consistent with
the ramification of which only one is defined over~$\Q$:
\be
t = -\frac{(256 x^3 + 224 x^2 + 232 x + 217)^2} {50000(x^2+1)},
\quad w_5(x) = \frac{24-7x}{7+24x}.
\label{14:5}
\ee
This yields the CM points of discriminants $-11$, $-35$, $-36$, $-84$
already known from $w_3$, and new points of discriminants $-91$,
$-100$, $-280$.

This accounts for eleven of the nineteen rational CM points
on~$\XX(1)$; the remaining ones were computed numerically as
we did for the $\Sigma=\{2,5\}$ curve.  We used the Schwarzian
equation
\be
t (16 t^2 + 13 t + 8) f''
+ (24 t^2 + 13 t + 4) f'
+ \left( \frac34 t + \frac{3}{16} \right) f = 0,
\label{schw.14}
\ee
for which the ``accessory parameter'' $3/16$ was again determined
by pulling back to $\XXo(3)$ and imposing the condition of
symmetry under~$w_3$.
We tabulate the coordinates $t=A/B$\/ and factorizations
for all nineteen points in Table~5.

\begin{figure}[t]
\centerline{\large Table 5}

\vspace*{2ex}

\centerline{
\def\M{\!\cdot\!}
\def\0{{\phantom0}}
\def\m{{^\0}}
$
\begin{array}{c|c|c|c|c|c}
|D|& D_0& D_1 &      A        &     B     &16 A^2 + 13 AB + 8 B^2\\
\hline
  4&  4 &  1  &      1        &     0     &         2^4        \\
  8&  8 &  1  &      0        &     1     &         2^3        \\
 11&  1 & 11  &     -1\0      &     1     &          11        \\
 35&  7 &  5  &     -5\0      &     1     &         7^3        \\
 36&  4 & 3^2 &  \0 4=2^2     &     3     &       2^2 11^2     \\
 84& 28 &  3  &  \0 4=2^2     & \0 9=3^2  &        2^2 7^3     \\
 72&  8 & 3^2 & \0 75=5^2 3   &   16=2^4  &       2^7 29^2     \\
 91&  7 & 13  &     -13       &   81=3^4  &       7^3 11^2     \\
 43&  1 & 43  &   -25=-5^2    &   81=3^4  &        29^2 43     \\
168& 56 &  3  &\m 200=2^3 5^2 & \0 9=3^2  &     2^4 7^3 11^2 \\
 88&  8 & 11  &\m-200=-2^3 5^2&   81=3^4  &      2^5 37^2 11   \\
100&  4 & 5^2 &\m-196=-2^2 7^2&  405=3^4 5&     2^2 11^2 43^2  \\
 67&  1 & 67  & -1225=-5^2 7^2&   81=3^4  &    11^2 53^2 67    \\
280& 56 &  5  &\m-845=-13^2 5 &1296=2^4 3^4&    2^8 7^3 11^2   \\
148&  4 & 37  &  1225=5^2 7^2 &5184=2^6 3^4&  2^4 11^2 67^2 37 \\
532& 28 & 19  &
   \0 96100=2^2 5^2 31^2 & 29241=3^4 19^2 & 2^2 7^3 11^2 29^2 37^2 \\
232&  8 & 29  &135200=2^5 5^2 13^2 & 194481=3^4 7^4 \0\0 &
		   2^3 11^2 53^2 109^2 29\\
427&  7 & 61  &-3368725=-5^2 47^2 61&6561=3^8 \0&7^3 11^2 29^2 43^2 53^2
\\
163&  1 & 163 &-2235025=-5^2 13^2 23^2 & 1185921=3^4 11^4 \0 &
                   37^2 107^2 149^2 163
\end{array}
$
}

\end{figure}

We see that $t$\/ is also a good coordinate 3-adically: a point
of $\XX(1)$ is supersingular at~3 iff the denominator of its
$t$-coordinate is a multiple of~3.  (It is supersingular at~5
iff $5|t$.)

\subsection{The case $\Sigma = \{3,5\}$}

Here the area of $\XX(1)$ is $1/3$.  This again is small enough to
show that there are only four elliptic points, but leaves two
possibilities for their indices: 2,2,2,6 or 2,2,3,3.  It turns
out that the first of these is correct.  This fact is contained
in the table of \cite[Ch.IV:2]{V}; it can also be checked as we
did in the cases $\Sigma=\{2,p\}$ ($p=3,5,7$) by exhibiting
appropriate elliptic elements of $\Gamma^*(1)$ --- which we need
to do anyway to compute the CM points.  We chose to write
write $\O=\Z[\frac12{1+c},e]$ with
\be
c^2+3 = e^2-5 = ce+ec = 0,
\label{cerel}
\ee
and found the elliptic elements
\be
s_2 = [4c-3e], s'_2 = [5c-3e-ce], s''_2 = [20c-9e-7ce], s_6 = [3+c]
\label{s2226}
\ee
[NB $20c-9e-7ce, 3+c\in 2\O$] of orders $2,2,2,6$ with
$s_2^{\phantom0} s'_2 s''_2 s_6^{\phantom0} = 1$.
The corresponding elliptic points
$P_2^{\phantom0},P'_2,P''_2,P_6^{\phantom0}$
have CM discriminants $-3,-12,-15,-60$.
For the first time we have a curve $\XXo(2)$, and here it
turns out that the elliptic points $P'_2$ is not ramified
in the cover $\XXo(2)/\XX(1)$: it admits two 2-isogenies
to itself, and one to $P''$.  Of the remaining elliptic points,
$P_6$ is totally ramified, and each of $P_2,P''_2$ has one
simple and one double preimage.  So we may choose coordinates
$x,t$\/ on $\XXo(2)$ and $\XX(1)$ such that $t=x(x-3)^2/4$,
with $t(P_6)=\infty$, $t(P_2)=0$, $t(P''_2)=1$.  To determine
$t(P'_2)$ we use the involution $w_2$, which switches $x=\infty$
(the triple pole) with $x=0$ (the simple zero), $x=4$ (the simple
preimage of $P''_2$) with one of the preimages~$x_1$ of $P'_2$ (the
one parametrizing the isogeny from $P'_2$ to $P''_2$), and
the other two preimages of $P'_2$ with each other.  Then
$w_2$ is $x\leftrightarrow 4x_1/x$, so the product of the roots of
$(t(x_1)-t(x))/(x-x_1)$ is $4x_1$.  Thus
\be
x(x-3)^2 - 4t(P'_2) = (x-x_1)(x^2+ax+4x_1)
\label{to15}
\ee
for some~$a$.  Equating $x^2$ coefficients yields $a=x_1-6$, and
equating the coefficients of~$x$ we find $9=10x_1^{\phantom0}-x_1^2$.
Thus $x_1=1$ or $x_1=9$; but the first would give us
$t(P'_2)=1=t(P''_2)$ which is impossible.  Thus $x_1=9$ and
$t(P'_2)=81$, with $w_2(x)=36/x$.  This lets us find six further
rational CM points, of discriminants $-7,-28,-40,-48,-120,-240$; we can
also solve for the accessory parameter $-1/2$ in the Schwarzian equation
\be
t (t-1) (t-81) f''
+ \left( \frac32 t^2 - 82 t + \frac{81}{2}\right) f'
+ \left( \frac1{18} t - \frac12 \right) f = 0 ,
\label{schw.15}
\ee
and use it to compute the remaining twelve rational CM points
numerically.  We tabulate the coordinates $t=A/B$\/ and
factorizations for the twenty-two rational CM points on $\XX(1)$
in Table~6.

\begin{figure}[t]
\centerline{\large Table 6}

\vspace*{2ex}

\centerline{
\def\M{\!\cdot\!}
\def\0{\phantom0}
\def\m{\phantom-}
$
\begin{array}{c|c|c|c|c|c|c}
|D|& D_0& D_1 &      A        &     B     &   A - B  & 81B - A \\ \hline
  3&  3 &  1  &      1        &     0     &     1    &   -1\m  \\
 12&  3 & 2^2 &      0        &     1     &   -1\m   &   3^4   \\
 60& 15 & 2^2 &      1        &     1     &     0    &  2^4 5  \\
 15& 15 &  1  &    81=3^4     &     1     &   2^4 5  &    0    \\
  7&  1 &  7  &   -27=-3^3    &     1     &  -2^2 7\m& 2^2 3^3 \\
 40&  5 & 2^3 &    27=3^3     &     2     &    5^2   &   3^3 5 \\
 43&  1 & 43  &   -27=-3^3    &   16=2^4  &    -43   &  3^3 7^2\\
195& 15 & 13  &    81=3^4     &   16=2^4  &   5\M13  &  3^3 5  \\
 48&  3 & 2^4 &   243=3^5\0   &     1     &   11^2 2 & -3^4 2\m\\
120& 15 & 2^3 &  -243=-3^5\0  &     2     &  -7^2 5\m&  3^4 5  \\
 28&  1 &2^2 7&  -675=-3^3 5^2&     1    &-2^2 13^2\m&2^2 3^3 7\\
115&  5 & 23  &   621=3^3 23  &   16=2^4  &   11^2 5 &  3^3 5^2\\
147&  3 & 7^2 &  -729=-3^6\0  & 112=2^4 7 &  -29^2\m &  3^4 11^2 \\
123&  3 & 41  &  2025=3^4 5^2 &   16=2^4  &   7^2 41 & -3^6\m    \\
 67&  1 & 67  &-3267=-3^3 11^2&   16=2^4  & -7^2 67\m&  3^3 13^2 \\
240& 15 & 2^4 &  9801=3^4 11^2&     1     &2^3 5^2 7^2&-2^3 3^5 5\m\\
267&  3 & 89  &  7225=5^2 17^2&   16=2^4  &  3^4 89  & -7^2 11^2\m \\
435& 15 & 29  & 21141=3^6 29\0&   16=2^4  & 5^3 13^2 &-3^4 5\M7^2\m\\
795& 15 & 53  & -6413=-11^2 53&432=2^4 3^3&-5\M37^2\m&5\M7^2 13^2\\
235&  5 & 47  &  1269=3^3 47\0&1024=2^{10}&  5\M7^2  &3^3 5^2 11^2\\
555& 15 & 37 &23409=3^4 17^2\0&1024=2^{10}&5\M11^2 37& 3^5 5\M7^2 \\
163&  1 &163 &-1728243=-3^3 11^2 23^2&1024=2^{10}&-103^2 163\m&
    3^3 7^2 37^2\\
\end{array}
$
}
\end{figure}

An equivalent coordinate that is also good 2-adically is
$(t-1)/4$, which is supersingular at~2 iff its denominator is even.

The elliptic curve $\Xone$ is obtained from $\XX(1)$ by extracting
square roots of $At$ and $B(t-1)(t-81)$ for some $A,B\in\Q^*/{\Q^*}^2$.
Using the condition that the Jacobian of $\Xone$, and any elliptic
curve occurring in the Jacobian of $\calXo(2)$, have conductor
at most 15 and 30 respectively, we find $A=B=-3$.  Then $\Xone$
has equation
\be
y^2 = -(3s^2+1) (s^2+27)
\label{15:X1}
\ee
(with $t=-3s^2$) and Jacobian isomorphic with elliptic curve
15C (15-A1); the curve intermediate between $\XX(2)$
and $\calXo(2)$ whose function field is obtained from $\Q(\XX(2))$
by adjoining $\sqrt{-3(t-1)(t-81)}$ has equation
\be
y^2=-3(x^4-10x^3+33x^2-360x+1296)
\label{15:X2}
\ee
and Jacobian 30C (30-A3).  Fundamental domains
for $\Gamma^*(1)$ and $\Gamma(1)$, computed by Michon~\cite{Mi}
and drawn by C.~L\'eger, can be found in \cite[pp.123--127]{V};
an equation for $\Xone$ birational with (\ref{15:X1}) is reported
in the table of~\cite[p.235]{JL}.

\subsection{The triangle group $G_{2,3,7}$ as an arithmetic group}
\def\0{^{\phantom0}}

It is well-known that the minimal quotient area of a discrete
subgroup of~$\Aut(\H)$ $=\PSL_2(\R)$ is $1/42$, attained only
by the triangle group $G_{2,3,7}$, and that the Riemann surfaces
$\H/\Gamma$ with $\Gamma$ a proper normal subgroup of finite index
in~$G_{2,3,7}$ are precisely the curves of genus $g>1$ whose number
of automorphisms attains Hurwitz's upper bound $84(g-1)$.
Shimura observed in~\cite{S2} that this group is arithmetic.\footnote{
  Actually this fact is due to Fricke \cite{F1,F2}, over a century ago;
  but Fricke could not relate $G_{2,3,7}$ to a quaternion algebra
  because the arithmetic of quaternion algebras had yet to be
  developed.
  }
Indeed, let $K$\/ be the totally real cubic field $\Q(\cos 2\pi/7)$
of minimal discriminant~$49$, and let $\A$ be a quaternion algebra
over~$K$\/ ramified at two of the three real places and at no finite
primes of~$K$.  Now for any totally real number field of degree~$n>1$
over~$\Q$, and any quaternion algebra over that field ramified at
$n-1$ of its real places, the group $\Gamma(1)$ of norm-1 elements of
a maximal order embeds as a discrete subgroup of
$\PSL_2(\R) = \Aut(\H)$, with $\H/\Gamma$ of finite area
given by Shimizu's formula
\be
{\rm Area}(\Xone) = \frac{d_K^{3/2}\zeta_K(2)}{4^{n-1}\pi^{2n}}
 \prod_{\wp\in\Sigma} (\Nm\wp-1)
\left[
 = \frac{(-1)^n}{2^{n-2}} \zeta_K(-1) \prod_{\wp\in\Sigma} (\Nm\wp-1)
\right]
\label{Shimizu}
\ee
(from which we obtained (\ref{Area1}) by taking $K=\Q$).  Thus,
in our case of $K=\Q(\cos 2\pi/7)$, $\Sigma=\{\infty,\infty'\}$,
the area of $\H/\Gamma(1)$ is $1/42$, so $\Gamma(1)$ must be isomorphic
with $G_{2,3,7}$.  From this Shimura deduced \cite[p.83]{S2} that for
any proper ideal $I\subset O_K$ his curve $\calX(I) = \H / \Gamma(I)$
attains the Hurwitz bound.  For instance, if $I$\/ is the prime
ideal $\wp\0_7$ above the totally ramified prime~$7$ of~$\Q$ then
$\calX(\wp\0_7)$ is the Klein curve of genus~3 with automorphism
group $\PSL_2(\F_7)$ of order 168.  The next-smallest example is
the ideal $\wp\0_8$ above the inert prime~2, which yields a curve of
genus~7 with automorphism group [P]SL$_2(\F_8)$ of order 504.
This curve is also described by Shimura as a ``known curve''$\!$,
and indeed it first appears in~\cite{FrickeM}; an equivalent curve
was studied in detail only a few years before Shimura by
Macbeath~\cite{Mac}, who does not cite Fricke, and the identification
of Macbeath's curve with Fricke's and with Shimura's $\calX(\wp\0_8)$
may first have been observed by Serre in a 24.vii.1990 letter to
Abhyankar.  At any rate, we obtain towers $\{\calX(\wp_7^r)\}_{r>0}$,
$\{\calX(\wp_8^r)\}_{r>0}$ of unramified abelian extensions which
are asymptotically optimal over the quadratic extensions of residue
fields\footnote{
  That is, over the fields of size $p^2$ for primes $p=7$ or
  $p\equiv\pm1\bmod7$, and $p^6$ for other primes~$p$.
  }
of~$K$\/ other than $\F_{49}$ and $\F_{64}$ respectively, which are
involved in the class field towers of exponents $7, 2$ of the Klein
and Macbeath curves over those fields.

These towers are the Galois closures of the covers of $\calX(1)$
by $\calXo(\wp_7^r)$, $\calXo(\wp_8^r)$, which again may be
obtained from the curves $\calXo(\wp\0_7)$, $\calXo(\wp\0_8)$ together
with their involutions.  It turns out that these curves both have
genus~0 (indeed the corresponding arithmetic subgroups
$\Gamma_0(\wp\0_7)$, $\Gamma_0(\wp\0_8)$ of $\Gamma(1)$ are the
triangle groups $G_{3,3,7}$, $G_{2,7,7}$ in \cite[class~X]{Takeuchi}).
The cover $\calXo(\wp\0_7)/\Xone$ has the same ramification data
as the degree-8 cover of classical modular curves
${\rm X}_0(7)/{\rm X}(1)$, and is thus given by the same
rational function
\bea
t &=& \frac
   {(x_7^4 - 8 x_7^3 - 18 x_7^2 - 88 x_7\0 + 1409)^2}
   {2^{13} 3^3 (9-x_7\0)}
\nonumber \\ \vspace*{-5.5ex} \label{Xp7} \\ \vspace*{-5.5ex}
  &=& 1 + \frac
   {(x_7^2 - 8 x_7\0 - 5)^3 (x_7^2 + 8 x_7\0 + 43)}
   {2^{13} 3^3 (9-x_7\0)}
\nonumber
\eea
(with the elliptic points of orders $2,3,7$ at $t=0,1,\infty$, i.e.\
$t$\/ corresponds to $1 - 12^{-3}j$).  The involution is different,
though: it still switches the two simple zeros $x\0_7=-4\pm\sqrt{-27}$
of $t-1$, but it takes the simple pole $x\0_7=0$ to itself instead of
the septuple pole at $x\0_7=\infty$.  Using (\ref{propA}) again we find
\be
w_{\wp\0_7}(x_7) = \frac{19 x_7\0 + 711} {13 x_7\0 - 19}.
\label{wp7}
\ee
For the degree-9 cover $\calXo(\wp\0_8)/\Xone$ we find
\bea
t &=& \frac
 {(1 - x_8) (2 x_8^4 + 4 x_8^3 + 18 x_8^2 + 14 x_8\0 + 25)^2}
 {27 (4 x_8^2 + 5 x_8\0 + 23)}
\nonumber \\ \vspace*{-5.5ex} \label{Xp8} \\ \vspace*{-5.5ex}
  &=& 1 - \frac
   {4(x_8^3 + x_8^2 + 5 x_8\0 - 1)^3}
   {27 (4 x_8^2 + 5 x_8\0 + 23)},
\nonumber
\eea
with the involution fixing the simple zero $x_8=1$ and switching
the simple poles, i.e.
\be
w_{\wp\0_8}(x_8) = \frac{51 - 19 x_8\0} {19 + 13 x_8\0}.
\label{wp8}
\ee
Note that all of these covers and involutions have rational
coefficients even though a priori they are only known to be
defined over~$K$.  This is possible because $K$\/ is a normal
extension of~$\Q$ and the primes $\wp\0_7$, $\wp\0_8$ used
to define our curves and maps are Galois-invariant.  To each
of the three real places of~$K$\/ corresponds a quaternion
algebra ramified only at the other two places, and thus
a Shimura curve $\Xone$ with three elliptic points $P_2,P_3,P_7$
to which we may assign coordinates $0,1,\infty$. Then Gal$(K/\Q)$
permutes these three curves; since we have chosen rational
coordinates for the three distinguished points, any point on or
cover of $\Xone$ defined by a Galois-invariant construction must
be fixed by this action of Galois and so be defined over~$\Q$.
The same applies to each of the triangle groups $G_{p,q,r}$
associated with quaternion algebras over number fields~$F$\/
properly containing~$\Q$, which can be found in cases III through XIX
of Takeuchi's list~\cite{Takeuchi}.  In each case, $F$\/ is Galois
over~$\Q$, and the finite ramified places of the quaternion algebra are
Galois-invariant.  Moreover, even when $G_{p,q,r}$ is not $\Gamma(1)$,
it is still related with $\Gamma(1)$ by a Galois-invariant construction
(such as intersection with $\Gamma_0(\wp)$ or adjoining $\ww_\wp$ or
$w_\wp$ for a Galois-invariant prime~$\wp$ of~$F$\/).  At least one of
the triangle groups in each commensurability class has distinct indices
$p,q,r$, whose corresponding elliptic points may be unambiguously
identified with $0,1,\infty$; this yields a model of the curve
$\H/G_{p,q,r}$, and thus of all its commensurable triangle curves,
that is defined over~$\Q$.

This discussion bears also on CM points on $\Xone$.  There are many
CM points on $\Xone$ rational over~$K$, but only seven of those
are $\Q$-rational: a CM point defined over~$\Q$ must come from a
CM field~$K'$ which is Galois not only over~$K$\/ but over~$\Q$.
Thus $K'$ is the compositum of $K$\/ with an imaginary quadratic field,
which must have unique factorization.  We check that of the nine such
fields only five retain unique factorization when composed with~$K$.
One of these, $\Q(\sqrt{-7}\,)$, yields the cyclotomic field
$\Q(e^{2\pi i/7})$, whose ring of integers is the CM ring for
the elliptic point $P_7: t=\infty$; two subrings still have
unique factorization and yield CM points $\wp\0_7$- and
$\wp\0_8$-isogenous to that elliptic point, which again are not
only $K\/$- but even $\Q$-rational thanks to the Galois invariance
of $\wp\0_7$, $\wp\0_8$.  The other four cases are the fields of
discriminant $-3,-4,-8,-11$, which yield one rational CM point each.
The first two are the elliptic points $P_3, P_2: t=1,0$.  To find the
coordinates of the CM point of discriminant~$-8$, and of the two points
isogenous with $P_7$, we may use the involutions (\ref{wp7},\ref{wp8})
on $\calXo(\wp\0_7)$ and $\calXo(\wp\0_8)$.  On $\calXo(\wp\0_7)$,
the involution takes $x\0_7=\infty$ to $19/13$, yielding the point
$t=3593763963/4015905088$ $\wp\0_7$-isogenous with $P_7$ on $\Xone$;
on $\calXo(\wp\0_8)$ the involution takes $x\0_8=\infty$ to $-19/13$,
yielding the point $t=47439942003/8031810176$
$\wp\0_8$-isogenous with $P_7$.  On the latter curve,
the second fixed point of the involution (besides $x\0_8=1$)
is $x\0_8 = -51/13$, which yields the CM point
$t=1092830632334/1694209959$ of discriminant~$-8$.
The two points isogenous with $P_7$ also arise from the
second fixed point of $w_{\wp\0_7}$ and a further solution
of $t(x\0_8) = t(w_{\wp\0_8}(x\0_8))$.

This still leaves the problem of locating the CM point
of discriminant $-11$.  We found it numerically using
quotients of hypergeometric functions as we did for $G_{2,4,6}$.
Let $c=2\cos 2\pi/7$, so $c$ is the unique positive root of
$c^3 + c^2 - 2c - 1$.  Consider the quaternion algebra over~$K$\/
generated by $i,j$ with
\be
i^2 = j^2 = c, \quad ij = -ji.
\label{ijrel}
\ee
This is ramified at the two other real place of~$K$, in which
$c$ maps to the negative reals $2\cos 4\pi/7$ and $2\cos 6\pi/7$,
but not at the place with $c=2\cos 2\pi/7$; since $c$ is a unit,
neither is this algebra ramified at any finite place with the
possible exception of~$\wp\0_8$, which we exclude using the fact
that the set of ramified places has even cardinality.  Thus
$K(i,j)$ is indeed our algebra~$\A$.  A maximal order~$\O$ is
obtained from $O_K[i,j]$ by adjoining the integral element
$(1 + ci + (c^2+c+1)j)/2$.  Then $\O^*$ contains the elements
$$
g_2 := ij/c, \qquad
g_3 := \frac12{(1 + (c^2-2) j + (3-c^2) ij)},
$$\vspace*{-4ex}\be\label{g237}\ee\vspace*{-4ex}$$
g_7 := \frac12(c^2+c-1 + (2-c^2) i + (c^2+c-2) ij)
$$
of norm~1, with $g_2^2 = g_3^3 = g_7^7 = -1$ and $g_2 = g_7 g_3$.
Thus the images of $g_2,g_3,g_7$ in $\Gamma(1)$ are elliptic elements
that generate that group.  A short search finds the linear
combination $(2-c^2)g_3 + (c^2+c) g_7 \in \O$ of discriminant $-11$;
computing its fixed point in~$\H$ and solving for~$t$\/ to high
precision (150 decimals, which turned out to be overkill), we obtain
a real number whose continued fraction matches that of
\be
\frac{88983265401189332631297917}{45974167834557869095293} =
\frac{7^3 43^2 127^2 139^2 207^2 659^2 11} {3^3 13^7 83^7}
\label{p11},
\ee
with numerator and denominator differing by
$2^9 29^3 41^3 167^3 281^3$.  Having also checked that this
number differs from the $t$-coordinates of the three non-elliptic
CM~points by products of small ($<10^4$) primes,\footnote{
  If $10^4$ does not seem small, remember that the factorizations
  are really over~$K$, not $\Q$; the largest {\em inert} prime
  that occurs is~19, and the split primes are really primes of~$K$\/
  of norm at most comparable with that of~19.
  }
and that it passes the supersingular test, we are quite confident
that (\ref{p11}) is in fact the $t$-coordinate of the CM point
of discriminant~$-11$.

\subsection{An irrational example: the algebras over
$\Q[\tau]/(\tau^3-4\tau+2)$ with $\Sigma=\{\infty_i,\infty_j\}$}

While our examples so far have all been defined over~$\Q$,
this is not generally the case for Shimura curves associated
to a quaternion algebra over a totally number field~$K$\/
properly containing~$\Q$.  For instance, $K$\/ may not be a
Galois extension of~$\Q$; or, $K$\/ may be Galois, but the
set of finite ramified places may fails to be Galois-stable;
or, even if that set is Galois-stable, the congruence conditions
on the subgroup of $\A^*/K^*$ may not be Galois-invariant, and
the resulting curve would not be defined over~$\Q$ even though
$\Xone$ would be.  In each case different real embeddings of
the field yield different arithmetic subgroups of $\PSL_2(\R)$
and thus different quotient curves.  We give here what is probably
the simplest example: a curve $\Xone$ associated to a quaternion
algebra with no finite ramified places over a totally real cubic
field which is not Galois over~$\Q$.  While the curve has genus~0,
no degree-1 rational function on it takes $\Q$-rational values at
all four of its elliptic points, and the towers of modular curves
over this $\Xone$ are defined over~$K$\/ but not over~$\Q$.

Let $K$\/ be the cubic field $\Q[\tau]/(\tau^3-4\tau+2)$ and
discriminant $148=2^2 37$, which is minimal for a totally real
non-Galois field.  Let $\A/K$\/ be a quaternion algebra ramified
at two of the three real places and at no finite primes of~$K$.
Using {\sc gp/pari} to compute $\zeta_K(2)$, we find from
Shimizu's formula~(\ref{Shimizu}) that the associated Shimura
curve $\Xone=\XX(1)$ has hyperbolic area $.16666\ldots$;
thus the area is $1/6$ and, since $\A$ is not in Takeuchi's
list, the curve $\Xone$ has genus~0 and four elliptic points, one of
order~3 and three of order~2.  The order-3 point $P_3$ has discriminant
$-3$ as expected, but the order-2 points are a bit more interesting:
their CM field is $K(i)$, but the ring of integers of that field
is not $O_K[i]$!  Note that the rational prime~2 is totally ramified
in~$K$, being the cube of the prime $(\tau)$; thus $(1+i)/\tau$
is an algebraic integer, and we readily check that it generates
the integers of $K(i)$ over $O_K$.  One of the elliptic points,
call it $P_2$, has CM ring $O_K[(1+i)/\tau]$ and discriminant
$-4/\tau^2$; of its three $(\tau)$-isogenous points, one is $P_2$
itself, and the others are the remaining elliptic points $P'_2,P''_2$,
with CM ring $O_K[i]$ of discriminant $-4$.

Thus the modular curve $\calXo((\tau))$ is a degree-3 cover of
$\Xone$ unramified above the elliptic point $P_2$, and ramified
above the other three elliptic points with type 3 for $P_3$ and
21 for $P'_2,P''_2$.  This determines the cover up to
$\bar K$-isomorphism --- the curve $\calXo((\tau))$ has genus~0,
and we can choose coordinates $x$ on that curve and $t$\/ on $\Xone$
such that $t(P_3)=\infty$ and $t = x^3 - 3 c x$ for some $c\neq0$ ---
but not the location of the unramified point $P_2$
relative to the other three elliptic points.  To
determine that we once again use the involution, this time
$w_{(\tau)}$, of $\calXo((\tau))$: this involution fixes the
point above $P_2$ corresponding to its self-isogeny, and
pairs the other two preimages of~$P_2$ with the simple
preimages of $P'_2,P''_2$.  We find that there are three
ways to satisfy this condition:
\be
t = x^3 - 3 (\tau^2-3) x, \quad
P_2: t = 1300 - 188 \tau - 351 \tau^2, \quad
P'_2, P''_2: t = \pm 2 (\tau^2-3)^{3/2},
\label{148:2}
\ee
and its Galois conjugates.  The correct choice is determined
by the condition that the Shimura curves must be fixed by
the involution of the Galois closure of $K/\Q$ that switches
the two real embeddings of~$K$\/ that ramify~$\A$: the image
of~$\tau$ under the the third (split) embedding must be used
in (\ref{148:2}).  We find that the simple and double preimages
of $P'_2,P''_2$ are $x=\pm 2 \sqrt{a^2-3}$, $\mp\sqrt{a^2-3}$,
and the preimages of $P_2$ are $12 - 2\tau - 3\tau^2$
(fixed by $w_{(\tau)}$) and
$(-12 + 2\tau + 3\tau^2 \pm (3a^2-12)\sqrt{a^2-3})/2$.
{}From this we recover as usual the tower of curves $\calXo((\tau)^r)$,
whose reductions at primes of~$K$\/ other than $\tau$ are
asymptotically optimal over the quadratic extensions of the primes'
residue fields, and which in this case is a tower of double (whence
cyclic) covers unramified above the genus-3 curve $\calXo((\tau)^4)$
and thus involved in that curve's class-field tower.

\subsection{Open problems}

\subsubsection{Computing modular curves and covers.}
Given a nonempty even set $\Sigma$ of rational primes, and thus
a quaternion algebra $\A/\Q$, how to compute the curve $\XX(1)$
together with its Schwarzian equation and modular covers such as
$\Xone$ and $\XXo(l)$?  Even in the simplest case $\Sigma=\{2,3\}$
where $\Gamma^*(1)$ is a triangle group and all the covers
$\XXo(l)/\XX(1)$ are in principle determined by their ramifications,
finding those covers seems at present a difficult problem once
$l$\/ gets much larger than the few primes we have dealt with here.
This is the case even when $l$\/ is still small enough that
$\XXo(l)$ has genus small enough, say $g\leqs5$, 
that the curve should have a simple model in projective space.
For instance, according to \ref{genus(6)} the curve $\XXo(73)$
has genus~1.  Thus its Jacobian is an elliptic curve; moreover
it must be one of the six elliptic curves of conductor $6\cdot 73$
tabulated in~\cite{Cr}.  Which one of those curves it is, and
which principal homogeneous space of that curve is isomorphic with
$\XXo(73)$, can probably be decided by local methods such as those
of~\cite{Ku}; indeed such a computation was made for $\calXo(11)$
in D.~Roberts' thesis~\cite{DaveRoberts}.
But that still leaves the problem of finding the degree-74 map
on that curve which realizes the modular cover $\XXo(73)\ra\XX(1)$.
For classical modular curves (i.e.\ with $\Sigma=\emptyset$)
of comparable and even somewhat higher levels, the equations
and covers can be obtained via $q$-expansions as explained
in~\cite{Modcomp}; but what can we do here in the absence of
cusps and thus of $q$-expansions?  Can we do anything at all
once the primes in $\Sigma$ are large or numerous enough to
even defeat the methods of the present paper for computing
$\XX(1)$ and the location of the elliptic points
on this curve?  Again this happens while the genus of $\XX(1)$ 
is still small; for instance it seems already a difficult problem
to locate the elliptic points on all curves $\XX(1)$ of genus~zero
and determine their Schwarzian equations, let alone find equations
for all curves $\XX(1)$ of genus 1, 2, or~3.  By \cite{I2} the
existence of the involutions $w_l$ on $\XXo(l)$ always suffices
in principle to answer these questions, but the computations
needed to actually do this become difficult very quickly;
it seems that a perspicuous way to handle these computations,
or a new and more efficient approach, is called for.  The reader
will note that so far we have said nothing about computing with
{\em modular forms}\/ on Shimura curves.  Not only is this an
intriguing question in its own right, but solving it may also allow
more efficient computation of Shimura curves and the natural maps
between them, as happens in the classical modular setting.

In another direction, we ask: is there a prescription, analogous
to (\ref{conds}), for towers of Shimura curves whose levels
are powers of a ramified prime of the algebra?  For a concrete
example (from case~III of~\cite{Takeuchi}), let $\A$ be the quaternion
algebra over $\Q(\sqrt2\,)$ with $\Sigma=\{\infty_1,\wp_2\}$,
where $\infty_1$ is one of the two Archimedean places and
$\wp_2$ is the prime ideal $(\sqrt2\,)$ above~$2$; let $\O\subset\A$
be a maximal order, $I=I_{\wp_2}\subset\O$ the ideal of elements
whose norm is a multiple of~$\sqrt2$, and
\be
\Gamma_n = \{ [a] \in O^*_1/\{\pm1\} : a \equiv 1 \bmod I^n \}
\label{Gamman}
\ee
for $n=0,1,2,\ldots$\ .
Then $\Gamma_{n+1}$ is a normal subgroup of $\Gamma_n$ with
index $3,2^2,2$ according as $n=0$, $n$ is odd, or $n$ is even
and positive.  Consulting \cite{Takeuchi}, we find that
$\Gamma_0,\Gamma_1$ are the triangle groups $G_{3,3,4}$ and
$G_{4,4,4}$.  Let $X_n$ be the Shimura curve $\H/\Gamma_n$,
which parametrizes principally polarized abelian fourfolds
with endomorphisms by~$\A$ and complete level-$I^n$ structure.
Then $\{X_n\}_{n>0}$ is a tower of $\Z/2$ or $(\Z/2)^2$
covers, unramified above the curve $X_3$.  Moreover, $X_n$ has
genus zero for $n=0,1,2$, while $X_3$ is isomorphic with the
curve $y^2=x^5-x$ of genus~2 with maximal automorphism group.
The reduction of this tower at any prime $\wp\neq\wp_2$ of
$\Q(\sqrt2\,)$ is asymptotically optimal over the quadratic
extension of the residue field of~$\wp$.  So we ask for explicit
recursive equations for the curves in this tower.  Note that unlike
the tower~(\ref{Xotower}), this one does not seem to offer a
$w_l$ or $\ww_{\wp_2}$ shortcut.

\subsubsection{CM points.}
Once we have found a Shimura modular curve together with a Schwarzian
equation, we have seen how to compute the coordinates of CM points
on the curve, at least as real or complex numbers to arbitrary
precision.  But this still leaves many theoretical and computational
questions open.  For instance, what form does the Gross-Zagier
formula~\cite{GZ} for the difference between $j$-invariants of elliptic
curves take in the context of Shimura curves such as $\XXo(1)$ or
$\Xone$?  Note that a factorization theorem would also yield a rigorous
proof that our tabulated rational coordinates of CM points are correct.
Our tables also suggest that at least for rational CM points the heights
increase more or less regularly with $D_1$; can this be explained
and generalized to CM points of degree $>1$?  For CM points on
the classical modular curve X$(1)$ this is easy: a CM $j$-invariant
is an algebraic integer, and its size depends on how close the
corresponding point of $\H/\PSL_2(\Z)$ is to the cusp; so for instance
if $\Q(\sqrt{-D})$ has class number~1 then the CM $j$-invariant of
discriminant $-D$\/ is a rational integer of absolute value
$\exp(\pi\sqrt{D})+O(1)$.  But such a simple explanation probably
cannot work for Shimura curves which have neither cusps nor
integrality of CM points.  Within a commensurability class
of Shimura curves (i.e.\ given the quaternion algebra~$\A$),
the height is inversely proportional to the area of the curve;
does this remain true in some sense when $\A$ is varied?

As a special case we might ask: how does the minimal polynomial
of a CM point of discriminant~$-D$\/ factor modulo the primes
contained in~$D_1$?  That the minimal polynomials for CM $j$-invariants
are almost squares modulo prime factors of the discriminant was a key
component of our results on supersingular reduction of elliptic curves
\cite{ss:Q,ss:K}; analogous results on Shimura curves may likewise yield
a proof that, for instance, for every $t\in\Q$ there are infinitely
many primes~$p$ such that the point on the $(2,4,6)$ curve with
coordinate~$t$ reduces to a supersingular point mod~$p$.

\subsubsection{Enumeration and arithmetic of covers.}
When an arithmetic subgroup of $\PSL_2(\R)$ is commensurable
with a triangle group $G=G_{p,q,r}$, as was the case for the
$\Sigma=\{2,3\}$ algebra, any modular cover $\H/G'$ of $\H/G$\/
(for $G'\subset G$\/ a congruence subgroup) is ramified above
only three points on the genus-0 curve $\H/G$.  We readily
obtain the ramification data, which leave only finitely many
possibilities for the cover.  We noted that, even when there
is only one such cover, actually finding it can be far from
straightforward; but much is known about covers of $\Pr^1$
ramified at three points --- for instance, the number of such
covers with given Galois group and ramification can be computed
by solving equations in the group (see \cite{Mat}), and the cover
is known~\cite{Beckmann} to have good reduction at each prime not
dividing the size of the group.
But when $G$, and any group commensurable with it, has positive
genus or more than three elliptic points, we were forced
to introduce additional information about the cover, namely
the existence of an involution exchanging certain preimages
of the branch points.  In the examples we gave here (and in
several others to be detailed in future work) this was enough
to uniquely determine the cover $\H/G' \ra \H/G$.  But there is
as yet no general theory that predicts the number of solutions
of this kind of covering problem.  The arithmetic of the solutions
is even more mysterious: recall for instance that in our final
example the cubic field $\Q[\tau]/(\tau^3-4\tau+2)$ emerged out
of conditions on the cover $\calXo((\tau))/\Xone$ in which that
field, and even its ramified prime~$37$, are nowhere to be seen.

\section{Appendix: Involutions of $\Pr^1$}

We collect some facts concerning involutions of the projective
line over a field of characteristic other than~2.  We do this
from a representation-theoretic point of view, in the spirit
of~\cite{FH}.  That is, we identify a pair of points $t_i=(x_i:y_i)$
$(i=1,2)$ of $\Pr^1$ with a binary quadric, i.e.\ a one-dimensional
space of homogeneous quadratic polynomials $Q(X,Y)=AX^2+2BXY+CY^2$,
namely the polynomials vanishing at the two points; we regard the
three-dimensional space $V_3$ of all such polynomials
$AX^2+2BXY+CY^2$ as a representation of the group SL$_2$
acting on $\Pr^1$ by unimodular linear transformations of $(X,Y)$.

An invertible linear transformation of a two-dimensional vector
space $V_2$ over any field yields an involution of the projective
line $\Pr^1=\Pr(V_2^*)$ if and only if it is
not proportional to the identity and its trace vanishes (the
first condition being necessary only in characteristic~2).
Over an algebraically closed field of characteristic other than~2,
every involution of $\Pr^1$ has two fixed points, and any two points
are equivalent under the action of $\PSL_2$ on $\Pr^1$.  It is
clear that the only involution fixing $0,\infty$ is $t\lra -t$;
it follows that any pair of points determines a unique involution
fixing those two points.  Explicitly, if $B^2\neq AC$, the involution 
fixing the distinct roots of $AX^2+2BXY+CY^2$ is
$(X:Y) \lra (BX+CY:-AX-BY)$.  Note that the 2-transitivity of $\PSL_2$
on $\Pr^1$ also means that this group acts transitively on the
complement in the projective plane $\Pr V_3$ of the conic $B^2=AC$\/ 
(and also acts transitively on that conic); indeed it is well-known
that $\PSL_2$ is just the special orthogonal group for the discriminant
quadric $B^2-AC$ on~$V_3$.

Now let $Q_1,Q_2\in V_3$ be two polynomials without a common zero.
Then there is a unique involution of $\Pr^1$ switching the roots
of $Q_1$ and also of $Q_2$.  (If $Q_i$ has a double zero the condition
on $Q_i$ means that its zero is a fixed point of the involution.)
This can be seen by using the automorphism group $\Aut(\Pr^1)=$PGL$_2$
to map $Q_i$ to $XY$ or $Y^2$ and noting that the involutions that
switch $t=0$ with $\infty$ are $t\lra a/t$ for nonzero~$a$, while the
involutions fixing $t=\infty$ are $t\lra a-t$ for arbitrary~$a$.
As before, we regard the involution determined in this way by $Q_1,Q_2$
as an element of $\Pr V_3$.  This yields an algebraic map~$f$ from
(an open set in) $\Pr V_3 \times \Pr V_3$, parametrizing $Q_1,Q_2$
without common zeros, to $\Pr V_3$.  We next determine this map
explicitly.

First we note that this map is covariant under the action of $\PSL_2$:
we have $f(gQ_1,gQ_2)=g(f(Q_1,Q_2))$ for any $g\in\PSL_2$.  Next we show
that $f$\/ has degree~1 in each factor.  Using the action of $\PSL_2$
it is enough to show that if $Q_1=XY$ or $Y^2$ then $f$\/ is linear
as a function of $Q_2=AX^2+2BXY+CY^2$.  In the first case, the
involution is $t\lra C/At$ and its fixed points are the roots of
$AX^2-CY^2$.  In the second case, the involution is $t\lra(-2B/A)-t$
with fixed points $t=\infty$ and $t=-B/A$, i.e.\ the roots of
$AXY+BY^2$.  In either case the coefficients of $f(Q_1,Q_2)$ are
indeed linear in $A,B,C$.

But it turns out that these two conditions completely determine~$f$\/:
there is up to scaling a unique $\PSL_2$-covariant bilinear map from
$V_3\times V_3$ to $V_3$; equivalently, $V_3$ occurs exactly once
in the representation $V_3\otimes V_3$ of~$\PSL_2$.  In fact it is
known (see e.g.\ \cite[\S 11.2]{FH})
that $V_3\otimes V_3$ decomposes as $V_1 \oplus V_3 \oplus V_5$,
where $V_1$ is the trivial representation and $V_5$ is the space of
homogeneous polynomials of degree~4 in $X,Y$.
The factor $V_3$ is particularly easy to see, because it is just
the antisymmetric part $\bigwedge^2 V_3$ of $V_3\otimes V_3$.
Now the next-to-highest exterior power $\bigwedge^{\!\dim V - 1} V$\/
of any finite-dimensional vector space~$V$\/ is canonically isomorphic
with $(\det V) \otimes V^*$, where $\det V$\/ is the top exterior power
$\bigwedge^{\!\dim V} V$.
Taking $V=V_3$, we see that $\det V_3$ is the trivial representation
of~$\PSL_2$.  Moreover, thanks to the invariant quadric $B^2-AC$\/
we know that $V_3$ is self-dual as a $\PSL_2$ representation.
Unwinding the resulting identification
$\bigwedge^2 V_3 \stackrel\sim\ra V_3^* \stackrel\sim\ra V_3$,
we find:

{\bf Proposition A}.  {\sl Let $Q_i = A_i X^2 + 2B_i XY + C_i Y^2$
($i=1,2$) be two polynomials in~$V_3$ without a common zero.
Then the unique involution of $\Pr^1\!$ switching the roots
of $Q_1$ and also of $Q_2$ is the involution whose fixed points
are the roots of
\be
(A_1 B_2 - A_2 B_1) X^2 + (A_1 C_2 - A_2 C_1) XY + (B_1 C_2 - B_2 C_1) Y^2,
\label{w}
\ee
i.e.\ the fractional linear transformation
\be
t \llra \frac
 {(A_1 C_2 - A_2 C_1)t + 2(B_1 C_2 - B_2 C_1)}
 {2(B_1 A_2 - B_2 A_1)t + (C_1 A_2 - C_2 A_1)}.
\label{propA}
\ee
}

{\sl Proof}\/: The coordinates of $Q_1 \wedge Q_2$ for the basis
of $V_3^*$ dual to $(X^2,2XY,Y^2)$ are ($B_1 C_2 - B_2 C_1$,
$A_2 C_1 - A_1 C_2$, $A_1 B_2 - A_2 B_1$).  To identify $V_3^*$
with $V_3$ we need a $\PSL_2$-invariant element of $V_3^{\otimes 2}$.
We could get this invariant from the invariant quadric
$B^2 - AC \in V_3^{*\otimes 2}$, but it is easy enough
to exhibit it directly: it is
\be
X^2 \otimes Y^2 - \frac12\, 2XY \otimes 2XY + Y^2 \otimes X^2,
\label{VVinv}
\ee
the generator of the kernel of the multiplication map
Sym$^2(V_3)\ra V_5$.  The resulting isomorphism from $V_3^*$
to~$V_3$ takes the dual basis of $(X^2,2XY,Y^2)$ to
$(Y^2,-XY,X^2)$, and thus takes $Q_1\wedge Q_2$ to (\ref{w})
as claimed.~~$\Box$

Of course this is not the only way to obtain (\ref{propA}).
A more ``geometrical'' approach (which ultimately amounts to
the same thing) is to regard $\Pr^1$ as a conic in $\Pr^2$.
Then involutions of~$\Pr^1$ correspond to points $p\in\Pr^2$
not on the conic: the involution associated with~$p$ takes
any point~$q$ of the conic to the second point of intersection
of the line $pq$ with the conic.  Of course the fixed points
are then the points~$q$ such that $pq$ is tangent to the conic
at~$q$.  Given $Q_1,Q_2$ we obtain for $i=1,2$ the secant of
the conic through the roots of $Q_i$, and then $p$ is the
intersection of those secants.

{}From either of the two approaches we readily deduce

{\bf Corollary B}.  {\sl Let $Q_i = A_i X^2 + 2B_i XY + C_i Y^2$
($i=1,2,3$) be three polynomials in~$V_3$ without a common zero.
Then there is an involution of $\Pr^1\!$ switching the roots of
$Q_i$ for each $i$ if and only if the determinant
\be
\left| \matrix{
A_1 & B_1 & C_1 \cr
A_2 & B_2 & C_2 \cr
A_3 & B_3 & C_3 \cr
} \right|
\label{CorB}
\ee
vanishes.}

As an additional check on the formula (\ref{w}), we may compute that
the discriminant of that quadratic polynomial is exactly the resolvent
\be
\det \left( \begin{array}{cccc}
A_1 & 2B_1 &  C_1 &  0  \\
 0  &  A_1 & 2B_1 & C_1 \\
A_2 & 2B_2 &  C_2 &  0  \\
 0  &  A_2 & 2B_2 & C_2
\end{array}
\right)
\label{resolvent}
\ee
of $Q_1,Q_2$ which vanishes if and only if these two polynomials
have a common zero.

\pagebreak

\centerline{\bf Corrigendum}

\vspace*{2ex}

While I was preparing my paper (call it [SCC]) I did not have David
Roberts' thesis [Ro] to hand.  Roberts has now kindly provided a copy
of [Ro]; it turns out that the second-hand information I had to rely on
concerning the contents of this thesis was wrong in several details,
requiring specific corrections as follows.

I ask several times for a formula for the factorization of the
differences between the coordinates of CM points on Shimura curves
analogous to the formulas proved Gross and Zagier for X$_0(1)$ [GZ]
and obtained experimentally by Yui and Zagier for singular values
of the Weber functions~[YZ].
Roberts had already answered this question in principle,
at least for the ${\cal X}_0$ curves, by obtaining an arithmetic
intersection formula for CM points on these curves (section 6.5
of [Ro]).  It would still take some work to extract from it (say)
an explicit factorization of the difference between two such points
on a Shimura curve of genus zero, but the technical framework now
exists.

Roberts identified arithmetically most of the elliptic curves of
conductor 60 or less which arise as Jacobians of Shimura curves
for quaternion algebras over~{\bf Q} [Ro, \S7.4].  This extends
considerably the list of previously computed Shimura curves,
and includes most that arise in [SCC].  Contrary to
a statement in the first paragraph of [SCC, \S5.5], the
list in Roberts' thesis does not, however, include the curve of
conductor 66 which arises as the Jacobian of ${\cal X}^*_0(11)$
[ $= X_{11,6}$ in Roberts' notation], though it would be easy
to obtain from his methods.

%In particular, the Jacobian of the curve ${\cal X}_0(5)$ for the
%$\{2,3\}$ algebra [Roberts' $X_{5,6}$] is given wrongly at the end
%of [SCC,\S3.2]: it is 30-F~(A6), not the \hbox{2-isogenous} curve
%30-H~(A8) as I wrongly wrote.  The latter curve is the Jacobian
%of the genus-1 curve intermediate between ${\cal X}^*_0(5)$ and
%${\cal X}_0(5)$; I neglected to adjoin a further square root
%to obtain ${\cal X}_0(5)$ itself.

In [SCC,\S5.2], the sentence preceding (76) needs some
explanation: conductor ``at most'' 15 and 30 rather than exactly?
In fact Jac(${\cal X}(1)$) is known {\em a~priori} to have conductor
exactly~15, but Jac(${\cal X}_0(2)$) has factors of conductor~15
[namely Jac(${\cal X}(1)$), from ``oldforms'' on ${\cal X}_0(2)$]
as well as 30; if there were no newforms at all on ${\cal X}_0(2)$
its Jacobian would consist entirely of curves of conductor 15.
[Cf.\ the case of the classical modular curve X$_0(22)$, whose
Jacobian is isogenous with the square of X$_0(11)$.]
For that matter, the fact that the conductors are thus bounded
at all needs a reference.  There are various ways of doing this;
surely the easiest (given existing work) is to cite [Ro] for
results finding (factors isogenous to) the Jacobians of Shimura
curves inside the Jacobians of classical modular curves X$_0(N)$,
and using the known results about elliptic curves occurring in $J_0(N)$.

Thanks again to David Roberts for bringing these matters to my
attention.

\end{document}